\documentclass[10pt]{article}
\usepackage[latin1]{inputenc}
\usepackage{epsfig}
\usepackage{color}
\usepackage[british,english]{babel}
\usepackage{amsthm}
\usepackage{amsmath}
\usepackage{amsfonts}
\usepackage{amssymb}
\usepackage{graphicx}
\newtheorem{corollary}{Corollary}[section]

\newtheorem{lemma}[corollary]{Lemma}
\newtheorem{proposition}[corollary]{Proposition}
\newtheorem{remark}[corollary]{Remark}
\newtheorem{theorem}[corollary]{Theorem}
\newfont{\sBlackboard}{msbm10 scaled 900}

\newcommand{\mylabel}[1]{\label{#1}
            \ifx\undefined\stillediting
            \else \fbox{$#1$}\fi }
\newcommand{\BE}{\begin{equation}}

\newcommand{\EEQ}{\end{equation}}
\newcommand{\rfb}[1]{\mbox{\rm
   (\ref{#1})}\ifx\undefined\stillediting\else:\fbox{$#1$}\fi}

\newfont{\Blackboard}{msbm10 scaled 1200}

\newfont{\roma}{cmr10 scaled 1200}

\def\CC{\rm \hbox{C\kern-.56em\raise.4ex
         \hbox{$\scriptscriptstyle |$}\kern+0.5 em }}

\def\n{|\kern -.05cm{|}\kern -.05cm{|}}

\def \noame{\noalign{\medskip}}
\newcommand{\mm}    {{\hbox{\hskip 0.5pt}}}

\newcommand{\bluff} {{\hbox{\raise 15pt \hbox{\mm}}}}

\newcommand{\ep}   {\varepsilon}
\usepackage{fancyhdr}

\makeatletter
\def\section{\@startsection {section}{1}{\z@}{-3.5ex plus -1ex minus
    -.2ex}{2.3ex plus .2ex}{\large\bf}}
\makeatother
\def\be{\begin{equation}}
\def\ee{\end{equation}}

\date{ }
\begin{document}
\thispagestyle{empty}
\title{\Large \bf On the effects of surface roughness in non-isothermal porous medium flow}\maketitle
\vspace{-2cm}
\begin{center}
Mar\'ia Anguiano\footnote{Departamento de An\'alisis Matem\'atico. Facultad de Matem\'aticas. Universidad de Sevilla. 41012-Sevilla (Spain) anguiano@us.es}, Igor Pa\v zanin\footnote{Department of Mathematics, Faculty of Science, University of Zagreb, Bijeni${\rm \check{c}}$ka 30, 10000 Zagreb (Croatia) pazanin@math.hr} and Francisco J. Su\'arez-Grau\footnote{Departamento de Ecuaciones Diferenciales y An\'alisis Num\'erico. Facultad de Matem\'aticas. Universidad de Sevilla. 41012-Sevilla (Spain) fjsgrau@us.es}
 \end{center}
\ \\
 \renewcommand{\abstractname} {\bf Abstract}
\begin{abstract}
We analyze a non-isothermal Darcy-Brinkman thin-film flow with a periodically oscillating boundary and viscous dissipation acting as a heat source. Using asymptotic analysis and the periodic unfolding method, we establish the convergence of velocity, pressure, and temperature fields as the small parameter (related to the film thickness and the period of the roughness) tends to zero. The limit problems depend on the relative scaling of the roughness wavelength and consist of coupled elliptic systems combining Reynolds-type equations with Darcy-Brinkman cell problems and reduced energy equation. In the critical roughness regime, the effective model exhibits a strong coupling induced by the oscillatory geometry, which does not occur in a smooth-boundary case.
\end{abstract}
\bigskip\noindent

\noindent {\small \bf AMS classification numbers:}. 35B27, 35Q35, 76S05.\\

\noindent {\small \bf Keywords:} Thin-film flow; surface roughness; viscous dissipation; Darcy-Brinkman system; unfolding method.
\ \\
\ \\
\section {Introduction}\label{S1}
Flows in thin fluid layers arise in a wide range of natural and technological applications, including lubrication systems, coating processes, porous and fractured media, microfluidic devices, and heat-transfer technologies (see e.g.~\cite{Szeri,Tabel}). In such configurations, one spatial dimension is significantly smaller than the others, which leads to strong anisotropy and motivates the derivation of reduced models via asymptotic analysis. Classical lubrication theory and Reynolds-type equations represent typical examples of such reductions, capturing the dominant macroscopic behavior of pressure and velocity fields in thin geometries. From a mathematical point of view, these effective models are rigorously obtained as limits of the Navier-Stokes equations or related viscous flow systems posed in thin domains (see e.g.~\cite{Bayada_Chambat,Mikelic2,PanasenkoPil}).\\
\ \\
In real-life situations, the confining surfaces of thin domains are rarely smooth and, thus, the effects of rough or rapidly oscillating boundaries on viscous flows has been extensively studied. In the context of thin-film settings, homogenization techniques are very effective, leading to modified Reynolds equations, effective permeability tensors, or additional source terms depending on the roughness scaling (see e.g.~\cite{Anguiano_SG,grau1,grau2,Pazanin_SG2}). A fundamental role in these analyses is played by modern homogenization tools, in particular the periodic unfolding method (introduced in \cite{Ciora,Ciora2,Cioran-book}), which allows one to rigorously treat oscillatory geometries and to identify different roughness regimes.\\
Thermal effects are often essential in thin-domain flows, since the resulting temperature variations can influence viscosity, permeability, and flow structure. Mathematical studies of non-isothermal thin-domain and porous-media flows include existence results, qualitative analyses, and asymptotic derivations of effective models, see e.g.~\cite{Boukrouche,Tesis_El_Mir,Applicable,CMAT,Pazanin_SG3}. Related investigations have also shown that viscous dissipation, as the irreversible conversion of a fluid's mechanical energy into thermal energy due to internal friction (viscosity) can significantly alter the temperature distribution. We refer the reader to e.g.~\cite{Chapter,Ciuperca,Multiscale}.\\
\ \\
Another important modeling aspect arises when the fluid interacts with permeable or porous structures, such as porous coatings or microstructured surfaces. In these situations, Darcy or Darcy-Brinkman type models provide a natural macroscopic description. The Darcy law \cite{Darcy}, later extended by Brinkman to include viscous effects \cite{Brinkman}, interpolates between Darcy and Stokes regimes and has been extensively studied both physically and mathematically, since the pioneering works \cite{Sanchez,Levy,Allaire}. In thin-domain settings, rigorous asymptotic analyses show that permeability effects may survive in the limit and lead to nontrivial reduced equations, allowing for a precise comparison between Darcy and Brinkman laws (see e.g.~\cite{AMC,PazaninRadulovic,Pazanin_SG2}). Recently, the non-isothermal thin-domain flows governed by Darcy-Brinkman type systems have also been addressed (see \cite{JEM,ZAMM,Pazanin_SG_Darcy}) providing new effective models governing the fluid flow. These works provide a solid mathematical foundation for non-isothermal Darcy-Brinkman models in thin geometries; however, they either focus on smooth domains or neglect the viscous dissipation effects.\\
\ \\
The purpose of the present work is to analyze the combined effects of surface roughness and viscous dissipation in a thin-domain flow governed by a non-isothermal Darcy-Brinkman system. In view of that, we consider a thin domain bounded by a flat lower surface and a periodically rough upper boundary. We impose a heat-flux type thermal boundary condition and aim to investigate the interaction between roughness and viscous dissipation under precise roughness-period scalings. Due to boundary roughness, we cannot use the approach from \cite{Pazanin_SG_Darcy} anymore, but we employ the periodic unfolding method adapted to thin domains with oscillatory boundaries. As a result, we rigorously derive effective limit models for the velocity, pressure, and temperature fields. More precisely, in the regime of relatively large roughness period, the limit system consists of Reynolds-type equations coupled with Darcy-Brinkman cell problems and an effective one-dimensional temperature equation driven by averaged viscous dissipation. In the critical roughness regime, there is a stronger coupling between the flow, temperature, and surface geometry, revealing roughness-induced effects that are absent in smooth-boundary thin-domain models. We believe that the obtained results may be used to enhance engineering models of porous medium flows involving rough surfaces, where thermal effects due to viscous dissipation play a significant role.\\
\ \\
The paper is organized as follows. In Section 2, we describe the thin domain with a periodically oscillating upper boundary, formulate the governing system, and establish uniform a priori estimates for weak solutions in appropriate Sobolev spaces. Section 3 is devoted to the construction of periodic unfolding operators adapted to the thin geometry, and to the derivation of uniform bounds and compactness properties for the unfolded sequences of the unknowns. In Section 4, we pass to the limit in the unfolded variational formulations of the momentum and energy equations, identify the limit fields in the physical variables and in the local periodic cells, and analyze their dependence on the roughness scaling. Section 5 completes the analysis by deriving, for each roughness regime under consideration, the corresponding macroscopic equations for the velocity, pressure and temperature together, explicitly characterizing the additional terms generated by viscous dissipation and boundary oscillations.
\noindent

\section{Preliminaries and problem description}\label{sec:setting}
The considered thin domain with an oscillating boundary is defined by
\begin{equation}\label{Omegaep}
\Omega^\ep=\{x=(x',x_3)\in\mathbb{R}^2\times \mathbb{R}\,:\, x'\in \omega,\ 0<x_3< h_\ep(x')\},
\end{equation}
where  $\omega$ is a open connected subset of $\mathbb{R}^2$ with Lipschitz boundary.  The function $h_\ep(x')= \varepsilon h\left(x'/\varepsilon^\ell\right)$ represents the real distance between the two surfaces. The small parameter $0<\varepsilon\ll 1$ is related to the film thickness and  $0<\varepsilon^\ell\ll 1$ is the period of the roughness. We primarily focus on the settings in which the period of the roughness is of order greater or equal than $\varepsilon$, i.e.~$0<\ell\leq 1$. In view of that, we analyze the following two cases\footnote{The highly oscillating case $\ell>1$ leads to the conclusion that, at macroscopic level, the problem reduces to the study of the fluid flow in a thin domain without roughness performed in \cite{Pazanin_SG_Darcy}. For reader's convenience, we recall these results in Section 5.3.}:
\begin{equation}\label{RelationAlpha}
\left\{\begin{array}{l}
\ep^\ell=\ep\quad\hbox{in the case  }\ell=1,\\
\noame
\displaystyle \ep\ll \ep^\ell \quad\hbox{in the case  }0<\ell<1\ \left(i.e.\quad \lim_{\ep\to 0}\ep^{1-\ell}=0\right).
\end{array}\right.
\end{equation}

\noindent Function $h$ is a strictly positive, bounded  and  $W^{1,\infty}$-function defined for $z'$ in $\mathbb{R}^2$, $Z'$-periodic with $Z'=(-1/2,1/2)^2$ the cell of periodicity in $\mathbb{R}^2$, and there exist $h_{\rm min}$ and $h_{\rm max}$ such that
$$0<h_{\rm min}=\min_{z'\in Z'} h(z'),\quad  h_{\rm max}=\max_{z'\in Z'}h(z')\,.$$

\noindent
The boundaries of $\Omega_\varepsilon$ are defined as follows
$$\begin{array}{c}
\displaystyle \Gamma_0 =\omega\times\{0\},\quad \Gamma_1^\varepsilon=\left\{(x',x_3)\in\mathbb{R}^2\times\mathbb{R}\,:\, x'\in \omega,\ x_3=  h_\varepsilon(x')\right\},\quad
\displaystyle\Gamma_{\rm lat}^\varepsilon=\partial\Omega_\ep\setminus (\Gamma_0\cup\Gamma_1^\ep).
\end{array}$$
For the microstructure of the periodicity of the boundary,  we assume that the domain $\omega$ is divided by a mesh of size $\ep^\ell$: for $k'\in\mathbb{Z}^2$, each cell $Z'_{k',\ep^\ell}=\ep^\ell k'+\ep^\ell Z'$. We define $\mathcal{T}_\ep=\{k'\in\mathbb{Z}\,:\, Z'_{k',\ep^\ell}\cap\omega\neq \emptyset\}$. In this setting, there exists an exact finite number of periodic sets $Z'_{k',\ep^\ell}$ such that $k'\in T_\ep$.   Also, we define $Z_{k',\ep^\ell}=Z'_{k',\ep^\ell}\times (0,h(z'))$ and  $Z=Z'\times (0,h(z'))$, which is the reference cell in $\mathbb{R}^3$. We define the boundaries $\widehat \Gamma_0=Z'\times \{0\}$, $\widehat \Gamma_1=Z'\times\{h(z')\}$.

\noindent Applying a dilatation in the vertical variable, i.e. $z_3=x_3/\ep$, we define the following rescaled sets
 \begin{equation}\label{domains_tilde}\begin{array}{c}
 \displaystyle \widetilde \Omega^\varepsilon=\{(x',z_3)\in\mathbb{R}^2\times \mathbb{R}\,:\, x'\in \omega,\ 0<z_3< h(x'/\ep^\ell)\},\\
 \noame
  \widetilde\Gamma_1^\varepsilon=\{(x',z_3)\in\mathbb{R}^2\times \mathbb{R}\,:\, x'\in \omega,\ z_3=h(x'/\ep^\ell)\}\quad \widetilde \Gamma_{\rm lat}^\varepsilon=\partial\widetilde \Omega_\ep\setminus (\widetilde \Gamma_0\cup\widetilde \Gamma_1^\ep).
  \end{array}
  \end{equation}
The quantity  $h_{\rm max}$ allows us to define the extended sets $\Omega=\omega\times  (0, h_{\rm max})$   and $\Gamma_1=\omega\times \{h_{\rm max}\}$.


\noindent We denote by $C$ a generic constant which can change from line to line. Moreover, $O_\ep$ denotes a  generic quantity, which can change from line to line, devoted to tend to zero when $\ep\to 0$.\\

\noindent Throughout the paper, we adopt the following notation. For a vector function ${\bf v} =( {\bf v} ', v_{3})$ defined in $\Omega^\varepsilon$, with ${\bf v}'=(v_1, v_2)$, we denote by
$\mathbb{D}:\mathbb{R}^3\to \mathbb{R}^3_{\rm sym}$  the symmetric part of the velocity gradient given by
$$\mathbb{D}[{\bf v}]={1\over 2}(D{\bf v}+(D{\bf v})^T)=\left(\begin{array}{ccc}
\partial_{x_1}v_1 &   {1\over 2}(\partial_{x_1}v_2 + \partial_{x_2}v_1) &   {1\over 2}(\partial_{x_3}v_1 + \partial_{x_1}v_3)\\
\noame
 {1\over 2}(\partial_{x_1}v_2 + \partial_{x_2}v_1) & \partial_{x_2}v_2 &   {1\over 2}(\partial_{x_3}v_2 + \partial_{x_2}v_3)\\
 \noame
 {1\over 2}(\partial_{x_3}v_1 + \partial_{x_1}v_3)&   {1\over 2}(\partial_{x_3}v_2 + \partial_{x_2}v_3)& \partial_{x_3}v_3
\end{array}\right).$$

\noindent For a vector function $\widetilde {\bf  v} =(\widetilde {\bf  v}', \widetilde  v_{3})$  a scalar function $\widetilde \varphi$, both defined in $\Omega$, and deduced from ${\bf v}$ and $\varphi$ after a dilatation in the vertical variable, the following operators are introduced:
 $$\begin{array}{c}
 \displaystyle \Delta_{\varepsilon}  \widetilde  {\bf v} =\Delta_{x'} \widetilde  {\bf v} + \varepsilon^{-2}\partial_{z_3}^2 \widetilde  {\bf v} ,\quad   \displaystyle\Delta_{ \varepsilon}\widetilde \varphi=\Delta_{x'}\widetilde \varphi+ \ep^{-2}\partial^2_{z_3}\widetilde \varphi,\\
  \noame
  (D_{\ep}\widetilde {\bf v})_{ij}=\partial_{x_j}\widetilde {v}_i\ \hbox{ for }\ i=1,2,3,\ j=1,2,\quad   (D_{\ep}\widetilde {\bf v})_{i3}=\ep^{-1}\partial_{z_3}\widetilde {v}_i \ \hbox{ for }\ i=1,2,3,\\
  \noame
  \nabla_{\ep}\widetilde\varphi=(\nabla_{x'}\widetilde \varphi, \ep^{-1}\partial_{z_3}\widetilde\varphi)^t,\quad {\rm div}_{\varepsilon}(\widetilde  {\bf   v})={\rm div}_{x'}(\widetilde {\bf v}')+\ep^{-1}{\partial_{z_3}}\widetilde  v_{3}.
  \end{array}$$
We also define $\mathbb{D}_{\ep}[\widetilde  {\bf v}]$ as
 $$\mathbb{D}_{\ep}[\widetilde  {\bf v}]=\mathbb{D}_{x'}[\widetilde  {\bf v}]+\ep^{-1}\partial_{z_3}[\widetilde  {\bf v}]=\left(\begin{array}{ccc}
\partial_{x_1}\widetilde  v_1 &   {1\over 2}(\partial_{x_1}\widetilde  v_2 + \partial_{x_2}\widetilde  v_1) &   {1\over 2} (\partial_{x_1}\widetilde  v_3+_\ep^{-1}\partial_{z_3}\widetilde  v_1)\\
\noame
 {1\over 2}(\partial_{x_1}\widetilde  v_2 + \partial_{x_2}\widetilde  v_1) & \partial_{x_2}\widetilde  v_2 &   {1\over 2} (\partial_{x_2}\widetilde  v_3+_\ep^{-1}\partial_{z_3}\widetilde  v_2)\\
 \noame
 {1\over 2} (\partial_{x_1}\widetilde  v_3+\ep^{-1}\partial_{z_3}\widetilde  v_1)&   {1\over 2} (\partial_{x_2}\widetilde  v_3+\ep^{-1}\partial_{z_3}\widetilde  v_2)& \ep^{-1}\partial_{z_3}\widetilde  v_3
\end{array}\right),$$
 where
\begin{equation}\label{def_der_sym_1}
  \mathbb{D}_{x'}[{\bf v}]=\left(\begin{array}{ccc}
\partial_{x_1}v_1 &   {1\over 2}(\partial_{x_1}v_2 + \partial_{x_2}v_1) &   {1\over 2} \partial_{x_1}v_3\\
\noame
 {1\over 2}(\partial_{x_1}v_2 + \partial_{x_2}v_1) & \partial_{x_2}v_2 &   {1\over 2} \partial_{x_2}v_3\\
 \noame
 {1\over 2} \partial_{x_1}v_3&   {1\over 2} \partial_{x_2}v_3& 0
\end{array}\right),
 \  \partial_{z_3}[{\bf v}]=\left(\begin{array}{ccc}
0 &   0&   {1\over 2} \partial_{z_3}v_1\\
\noame
 0& 0 &   {1\over 2} \partial_{z_3}v_2\\
 \noame
 {1\over 2} \partial_{z_3}v_1&   {1\over 2} \partial_{z_3}v_2& \partial_{z_3}v_3
\end{array}\right).
\end{equation}

\noindent Finally, we introduce the following operators applied to ${\bf v}'$:
\begin{equation}\label{def_der_sym_2}
  \mathbb{D}_{x'}[{\bf v}']=\left(\begin{array}{ccc}
\partial_{x_1}v_1 &   {1\over 2}(\partial_{x_1}v_2 + \partial_{x_2}v_1) &   0\\
\noame
 {1\over 2}(\partial_{x_1}v_2 + \partial_{x_2}v_1) & \partial_{x_2}v_2 &  0\\
 \noame
0&   0& 0
\end{array}\right),
 \quad \partial_{z_3}[{\bf v}']=\left(\begin{array}{ccc}
0 &   0&   {1\over 2} \partial_{z_3}v_1\\
\noame
 0& 0 &   {1\over 2} \partial_{z_3}v_2\\
 \noame
 {1\over 2} \partial_{z_3}v_1&   {1\over 2} \partial_{z_3}v_2& 0
\end{array}\right).
\end{equation}


We assume that the domain $\Omega^\ep$ is filled with a fluid-saturated sparsely packed porous medium, thus the flow is governed by the Darcy-Brinkman system coupled with the heat equation as follows:
\begin{equation}\label{system_1}
\left\{\begin{array}{rl}
\displaystyle -2\mu_{\rm eff}\,{\rm div}(\mathbb{D}[{\bf u}_\ep])+{\mu\over K_\ep}{\bf u}_\ep+\nabla p_\ep={\bf f}_\ep&\hbox{ in }\Omega^\ep,\\
\noame
{\rm div}( {\bf u}_\ep)=0&\hbox{ in }\Omega^\ep, \\
\noame
\displaystyle-k \Delta T_\ep={\mu\over K_\ep}|{\bf u}_\ep|^2+2\mu_{\rm eff}|\mathbb{D}[{\bf u}_\ep]|^2 &\hbox{ in }\Omega^\ep.
\end{array}\right.
\end{equation}
The unknowns in the above system are: ${\bf u}_\ep$ (the filter velocity), $p_\ep$ (the pressure) and $T_\ep$ (the temperature). ${\bf f}_\ep$ is the momentum source term given in advance, $\mu$ is the dynamic viscosity coefficient, $ \mu_{\rm eff}$ stands for the effective viscosity of the Brinkman term, $K_\ep$ is the permeability of the porous medium, whereas $k$ is the thermal conductivity.
The viscous dissipation function $\Phi({\bf u}_\ep, \mu, \mu_{\rm eff}, K_\ep)={\mu\over K_\ep}|{\bf u}|^2+2\mu_{\rm eff}|\mathbb{D}[{\bf u}]|^2$ appearing in the equation (\ref{system_1})$_3$ was proposed in \cite{alhadhrami}, guaranties the correct asymptotic behavior for a wide range of permeability values $K_\ep$.
We remark that the superscript $\ep$ is added to stress the dependence of the solution and given parameters on the small parameter.\\
\ \\
We endow the system (\ref{system_1}) with the following boundary conditions:
\begin{equation}\label{BCBot1}
\begin{array}{l}
\displaystyle {\bf u}^\varepsilon=0,\quad T_\ep=0\quad \hbox{on}\ \Gamma_1^\ep\cup \Gamma_{\rm lat}^\ep,
\end{array}
\end{equation}
\begin{equation}\label{BCBot}
\begin{array}{l}
\displaystyle {\bf u}^\varepsilon=0,\quad  k{\partial T_\ep\over \partial n}=b_{\ep}\quad \hbox{on}\ \Gamma_0.
\end{array}
\end{equation}
As emphasized in the Introduction, we allow the heat flux across the bottom of the fluid domain (see (\ref{BCBot}))$_2$.\\
\ \\
\noindent We use the following scalings on the given data with respect to the small parameter $\ep$:
\begin{itemize}
\item[--] Following \cite{PazaninRadulovic,Pazanin_SG2}, we assume:
\begin{equation}\label{Kvalue}K_\ep=\ep^2 K,\quad\hbox{with}\quad K=\mathcal{O}(1).
\end{equation}
\item[--]
The external source function is $x_3$-independent and scaled as (see \cite{Boukrouche,Pazanin_SG3}):
\begin{equation}\label{externalforces}
{\bf f}_\varepsilon= \ep^{-2}({\bf f}'(x'),0),\quad \hbox{with}\quad f'\in L^2(\omega)^2.
\end{equation}
\item[--] The given function $b^{\ep}$ is scaled as (see \cite{Boukrouche, Pazanin_SG_Darcy}):
\begin{eqnarray}
  b_\ep=\ep^{-1}b&\hbox{with}& b=\mathcal{O}(1). \label{boundary_b}
 \end{eqnarray}
\end{itemize}
We work in the following functional framework on $\Omega^\ep$, where $1<r<+\infty$ with $1/r+1/r'=1$:
$$\begin{array}{l}
W^{1,r}_{\Gamma_1^\ep\cup\Gamma_{\rm lat}^\ep}(\Omega^\ep)=\{\psi\in W^{1,r}(\Omega^\ep)\,:\, \psi=0\hbox{ on }\Gamma_1^\ep\cup\Gamma_{\rm lat}^\ep\},\quad
\displaystyle L^{r}_0(\Omega^\ep)=\left\{\psi\in L^{r'}(\Omega^\ep)\,:\, \int_{\Omega^\ep}\psi\,dx=0\right\}.
\end{array}$$
{\bf Weak formulation.} By multiplying (\ref{system_1})$_1$ by ${\bf v}\in H^1_0(\Omega^\ep)^3$,  (\ref{system_1})$_2$ by $\psi \in H^1_0(\Omega^\ep)$ and (\ref{system_1})$_3$ by $\zeta\in W^{1,q'}_{ \Gamma_1^\ep\cup\Gamma_{\rm lat}^\ep}$, $q'=q/(q-1)$ with $q\in (1,3/2)$, respectively, after integration over $\Omega^\ep$ we obtain the variational form associated to (\ref{system_1})-(\ref{BCBot}):\\
\ \\
Find ${\bf u}_\ep\in H^1_0(\Omega^\ep)^3$, $p_\ep\in L^2_0(\Omega^\ep)$ and $T_\ep\in W^{1,q}_{\Gamma_1^\ep\cup\Gamma_{\rm lat}^\ep}(\Omega^\ep)$, with $q\in(1,3/2)$, such that
\begin{equation}\label{formvar}
\left\{\begin{array}{l}
\displaystyle 2\mu_{\rm eff}\int_{\Omega^\ep}\mathbb{D}[{\bf u}_\ep]:\mathbb{D}[{\bf v}]\,dx+{\mu\over  K}\ep^{-2}\int_{\Omega^\ep}{\bf u}_\ep\cdot {\bf v}\,dx-\int_{\Omega^\ep}p_\ep\,{\rm div}({\bf v})\,dx=\ep^{-2}\int_{\Omega^\ep}{\bf f}'\cdot {\bf v}'\,dx,
\\
\noame
\displaystyle \int_{\Omega^\ep}{\bf u}_\ep\cdot \nabla \psi\,dx=0,
\\
\noame
\displaystyle  k \int_{\Omega^\ep}\nabla T_\ep\cdot \nabla\zeta\,dx={\mu\over  K}\ep^{-2}\int_{\Omega^\ep} |{\bf u}_\ep|^2\zeta\,dx+2\mu_{\rm eff}  \int_{\Omega^\ep}|\mathbb{D}[{\bf u}_\ep]|^2\zeta\,dx+\int_{\Gamma_0}\ep^{-1}b\,\zeta\,d\sigma,
\end{array}\right.
\end{equation}
\noindent for every ${\bf v}\in H^1_{0}(\Omega^\ep)^3$,  $\psi\in H^1_{0}(\Omega^\ep)$ and $\zeta\in W^{1,q'}_{\Gamma_1^\ep\cup\Gamma_{\rm lat}^\ep}(\Omega^\ep)$.\\
\ \\
\begin{remark}We observe that $q'=1/(q-1)>3$, following \cite{Gallouet}, by Sobolev inequalities, $\zeta\in L^\infty(\Omega^\varepsilon)$ and the right-hand side of (\ref{formvar})$_3$ make sense.\\
The well-posedness can be deduced by directly adapting the proof from \cite[Theorem 1]{Boukrouche} and \cite[Chapitre 2, Th\'eor\`eme 2.2]{Tesis_El_Mir} (see also \cite[Theorem 2.4]{Ciuperca} ). In view of that, the system (\ref{system_1})-(\ref{BCBot}) admits at least one solution $({\bf u}_\ep, p_\ep, T_\ep)\in H^1_{0}(\Omega^\ep)^3\times L^2_0(\Omega^\ep)\times W^{1,q}_{\Gamma_1^\ep\cup\Gamma_{\rm lat}^\ep}(\Omega^\ep)$ with $q\in (1,3/2)$.
\end{remark}
 \noindent The objective of the present paper is to  derive  the limit model describing the asymptotic behavior of the process governed by (\ref{system_1})-(\ref{boundary_b}) depending on the value of $\ell$. We start by applying the dilatation in the variable $x_3$
\begin{equation}\label{dilatacion}
 z_3={x_3\over \ep}
\end{equation}
to get the functions defined in $\widetilde \Omega^\ep$. As a consequence, the system (\ref{system_1}) is rewritten as:
\begin{equation}\label{system_1_dil}
\left\{\begin{array}{rl}
\displaystyle -2\mu_{\rm eff}\,{\rm div}_{\ep}(\mathbb{D}_{\ep}[\widetilde {\bf u}_\ep])+{\mu\over K}\ep^{-2}\widetilde {\bf u}_\ep+\nabla_{\ep} \widetilde p_\ep={\bf f}_\ep  &\hbox{ in } \widetilde \Omega^\ep,\\
\noame
{\rm div}_{\ep}( \widetilde {\bf u}_\ep)=0&\hbox{ in } \widetilde \Omega^\ep, \\
\noame
\displaystyle-  k \Delta_{\ep}\widetilde  T_\ep ={\mu\over  K}\ep^{-2}|\widetilde {\bf u}_\ep|^2+2\mu_{\rm eff}|\mathbb{D}_{\ep}[\widetilde {\bf u}_\ep]|^2 &\hbox{ in } \widetilde \Omega^\ep,
\end{array}\right.
\end{equation}
with
\begin{equation}\label{BCBot1}
\begin{array}{l}
\displaystyle \widetilde{\bf u}^\varepsilon=0,\quad \widetilde T_\ep=0\quad \hbox{on}\   \widetilde \Gamma^\ep_1 \cup  \widetilde \Gamma_{\rm lat}^\varepsilon,
\end{array}
\end{equation}
\begin{equation}\label{BC_dil}
\begin{array}{l}
\displaystyle  \widetilde{\bf u}^\varepsilon=0,\quad    k \nabla_{\ep}\widetilde T_\ep\cdot n=\varepsilon^{-1}  b\quad \hbox{on}\ \Gamma_0.
\end{array}
\end{equation}

\noindent In the above system, the unknown functions are given by ${\bf \widetilde  u}^\varepsilon(x',z_3)={\bf u}^\varepsilon(x', \varepsilon z_3)$, $\widetilde  p^\varepsilon(x',z_3)=p^\varepsilon(x', \varepsilon z_3)$ and $\widetilde  T_\ep(x',z_3)=T_\ep(x', \varepsilon z_3)$ for a.e. $(x',z_3)\in  \widetilde \Omega^\ep$. Let us write the weak variational formulation of (\ref{system_1_dil})-(\ref{BC_dil}): \\

Find $\widetilde  {\bf u}_\ep\in H^1_{0}(\widetilde \Omega^\ep)^3$, $\widetilde  p_\ep\in L^2_0(\widetilde \Omega^\ep)$ and $\widetilde  T_\ep\in  W^{1,r}_{ \widetilde \Gamma^\ep_1\cup \widetilde \Gamma^\ep_{\rm lat}}(\widetilde \Omega^\ep)$ with $q\in (1,3/2)$, such that
\begin{equation}\label{formvar_dil}
\left\{\begin{array}{l}
\displaystyle 2\mu_{\rm eff}\int_{\widetilde \Omega^\ep}\!\!\mathbb{D}_{\ep}[\widetilde  {\bf u}_\ep]:\mathbb{D}_{\ep}[\widetilde  {\bf v}]\,dx'dz_3\!+\!{\mu\over  K}\ep^{-2}\int_{\widetilde \Omega^\ep}\!\!\widetilde  {\bf u}_\ep\cdot \widetilde  {\bf v}\,dx'z_3\!-\!\int_{\widetilde \Omega^\ep}\!\!\widetilde  p_\ep\,{\rm div}_{\ep}(\widetilde  {\bf v})\,dx'dz_3=\int_{\widetilde \Omega^\ep}\!\!\!\ep^{-2}{\bf f}'\cdot \widetilde  {\bf v}'\,dx'dz_3,\\
\noame
\displaystyle \int_{\widetilde \Omega^\ep}\widetilde  {\bf u}_\ep\cdot \nabla_{\ep}\widetilde  \psi\,dx'dz_3=0,
\\
\noame
\displaystyle \ep^2\int_{\widetilde \Omega^\ep}  k\nabla_{\ep}\widetilde  T_\ep\cdot \nabla_{\ep}\widetilde  \zeta\,dx'dz_3 ={\mu\over K} \int_{\widetilde \Omega^\ep} |\widetilde  {\bf u}_\ep|^2\widetilde  \zeta\,dx'dz_3+2\mu_{\rm eff}\ep^{2}\int_{\widetilde \Omega^\ep}|\mathbb{D}_{\ep}[\widetilde  {\bf u}_\ep]|^2\widetilde  \zeta\,dx'dz_3+\int_{\omega}   b\,\widetilde \zeta\,dx',\end{array}\right.
\end{equation}
\indent where $\widetilde  {\bf v}, \widetilde \psi$ and $\widetilde  \zeta$ are deduced from $ {\bf v},  \psi$ and $  \zeta$ given in (\ref{formvar}) via the change of variable (\ref{dilatacion}).


Hereinafter, we use the following spaces:
$$\begin{array}{l}
\displaystyle   W^{1,r}_{ \widetilde \Gamma^\ep_1\cup \widetilde \Gamma^\ep_{\rm lat}}(\widetilde \Omega^\ep)=\{\widetilde \psi\in W^{1,r}(\widetilde \Omega^\ep)\,:\, \psi=0\hbox{ on }  \widetilde \Gamma^\ep_1 \cup  \widetilde \Gamma_{\rm lat}^\ep\},\quad
\displaystyle L^{r}_0(\widetilde \Omega^\ep)=\left\{\psi\in L^{r'}(\widetilde \Omega^\ep)\,:\, \int_{\widetilde \Omega^\ep}\widetilde  \psi\,dx'dz_3=0\right\}.
\end{array}$$

\noindent  The goal is to describe the asymptotic behavior of this new sequences ${\bf \widetilde  u}_\ep$,  $\widetilde  p_\ep$ and $\widetilde  T_\ep$, as $\ep$  tends to zero. For this, we first recall the corresponding {\it a priori} estimates for velocity, temperature and pressure, whose derivation can be found in \cite[Propositions 1 and 2 and Corollary 1]{Pazanin_SG_Darcy}.

\begin{proposition}[Proposition 1 in \cite{Pazanin_SG_Darcy}]\label{lemma_estimates} Assume $q\in(1,3/2)$ and let $(\widetilde {\bf u}^\varepsilon, \widetilde T_\ep)$ be a solution of the dilated problem (\ref{system_1_dil})-(\ref{BC_dil}). Then, there hold the following estimates:
\begin{eqnarray}
\displaystyle
\|\widetilde {\bf  u}^\varepsilon\|_{L^2(\widetilde\Omega^\ep)^3}\leq C , &\displaystyle
\|D_{\varepsilon} \widetilde {\bf  u}^\varepsilon\|_{L^2(\widetilde\Omega^\ep)^{3\times 3}}\leq C\varepsilon^{-1},& \|\mathbb{D}_{\ep}[\widetilde {\bf u}_\varepsilon]\|_{L^2(\widetilde\Omega^\ep)^{3\times 3}}\leq C\varepsilon^{-1}, \label{estim_sol_dil1}
\\
\noame
\displaystyle
\|\widetilde  T_\ep\|_{L^q(\widetilde\Omega^\ep)}\leq C, &\displaystyle
\|\nabla_{\varepsilon} \widetilde  T_\ep\|_{L^q(\widetilde\Omega^\ep)^{3}}\leq C\varepsilon^{-1}.&  \label{estim_sol_T_dil1}
\end{eqnarray}
\end{proposition}
\begin{remark}[Extension of $\widetilde{\bf u}_\ep$ and $\widetilde {T}_\ep$ to $\Omega=\omega\times  (0, h_{\rm max})$] From the boundary conditions satisfied by ${\bf \widetilde u}^\varepsilon$ and $\widetilde T_\ep$ on the top boundary $\widetilde \Gamma_1^\ep$,  we can extend it  by zero in $\Omega\setminus \widetilde{\Omega}^{\varepsilon}$ and we denote the extension by  the same symbol. As consequence, the extended velocity and temperature satisfy the same estimates given in Proposition \ref{lemma_estimates}.
\end{remark}
\begin{proposition}[Proposition 2 and Corollary 1 in \cite{Pazanin_SG_Darcy}]\label{prop_decomposition}
The following decomposition for $p_\varepsilon\in L^2_0(\Omega^\varepsilon)$ holds
\begin{equation}\label{decompositionp}
p_\varepsilon=p_\varepsilon^0+p_\varepsilon^1,
\end{equation}
where $p_\varepsilon^0\in H^1(\omega)$, which is independent of $x_3$, and $p_\varepsilon^1\in L^2(\Omega^\varepsilon)$. Moreover, the following estimates hold
\begin{equation}\label{estim_p0p1}
\begin{array}{c}
\displaystyle
\|p_\varepsilon^0\|_{H^1(\omega)}\leq C\ep^{-2} ,\\
\\
\displaystyle   \|p_\varepsilon^1\|_{L^2(\Omega^\varepsilon)}\leq C\ep^{-{1\over 2}} ,\quad  \|\widetilde  p_\varepsilon^1\|_{L^2(\widetilde\Omega^\ep)}\leq C\ep^{-1},
\end{array}
\end{equation}
where  $\widetilde  p_\varepsilon^1$ is the rescaled function associated with $p_\varepsilon^1$ defined by $\widetilde  p_\ep^1(x',z_3)=p_\ep^1(x',\ep z_3)$ for a.e. $(x',z_3)\in \widetilde \Omega^\ep$.
\end{proposition}
\begin{proof} As explained in \cite[Proposition 2]{Pazanin_SG_Darcy}, the result is based on a decomposition of the pressure $p_\ep$ in two pressures $p_\ep^0\in H^1(\omega)$ and $p_\ep^1\in L^2(\Omega^\ep)$ given in \cite[Corollary 3.4 with $q=2$]{CLS2}, which gives the following estimate
$$\|p_\ep^0\|_{H^1(\omega)}\leq C^{-{3\over 2}}\|\nabla p_\ep\|_{H^{-1}(\Omega^\ep)^3},\quad \|p_\ep^1\|_{L^2(\Omega^\ep)}\leq C\|\nabla p_\ep\|_{H^{-1}(\Omega^\ep)^3}.$$
We remark that the result in \cite{CLS2} can be applied to $\Omega^\ep$ because $\omega\subset \mathbb{R}^2$ has a Lipschitz boundary and the function $h\in W^{1,\infty}$. Then, it is proved in \cite[Corollary 1]{Pazanin_SG_Darcy} that $\|\nabla p_\ep\|_{H^{-1}(\Omega^\ep)^3}\leq C\ep^{3\over 2}$, which finishes the proof.

\end{proof}

\section{Adaptation of the unfolding method}\label{sec:unfolding}
 The change of variable (\ref{dilatacion}) does not capture the oscillations of the domain $\widetilde\Omega^\ep$. To capture them, we use an adaptation of the unfolding method (see for instance \cite{Ciora, Ciora2}) introduced to this context in \cite{Anguiano_SG} (see also \cite{Anguiano_SG3}). In view of that, given the unknowns $\widetilde{\bf u}_{\varepsilon},  \widetilde T_\ep, p_\ep^0$ and $\widetilde p_\ep^1$, we define
\begin{equation}\label{phihat}
\begin{array}{l}
\displaystyle\widehat{\bf u}_{\varepsilon}(x^{\prime},z)=\widetilde{\bf u}_{\varepsilon}\left( {\varepsilon^\ell}\kappa\left(\frac{x^{\prime}}{{\varepsilon^\ell}} \right)+{\varepsilon^\ell}z^{\prime},z_3 \right),\quad \hbox{a.e. }(x',z)\in \omega\times Z,\\
\displaystyle\widehat T_{\varepsilon}(x^{\prime},z)=\widetilde T_{\varepsilon}\left( {\varepsilon^\ell}\kappa\left(\frac{x^{\prime}}{{\varepsilon^\ell}} \right)+{\varepsilon^\ell}z^{\prime},z_3 \right),\quad \hbox{a.e. }(x',z)\in \omega\times Z,\\
\noame
\displaystyle\widehat p_{\varepsilon}^0(x^{\prime},z')=p_\ep^0\left( {\varepsilon^\ell}\kappa\left(\frac{x^{\prime}}{{\varepsilon^\ell}} \right)+{\varepsilon^\ell}z^{\prime} \right),\quad \hbox{a.e. }(x',z')\in \omega\times Z',\\
\noame
\displaystyle
\widehat p_{\varepsilon}^1(x^{\prime},z)=\widetilde p_{\varepsilon}^1\left( {\varepsilon^\ell}\kappa\left(\frac{x^{\prime}}{{\varepsilon^\ell}} \right)+{\varepsilon^\ell}z^{\prime},z_3 \right),\quad \hbox{ a.e. }(x^{\prime},z)\in \omega\times Z,
\end{array}\end{equation}
assuming the unknowns  be extended by zero outside $\omega$, where the function $\kappa:\mathbb{R}^2\to \mathbb{Z}^2$ is defined by
$$\kappa(x')=k'\Longleftrightarrow x'\in Z'_{k',1},\quad\forall\,k'\in\mathbb{Z}^2.$$

\begin{remark}\label{remarkCV}We make the following comments:
\begin{itemize}
\item The function $\kappa$ is well defined up to a set of zero measure in $\mathbb{R}^2$ (the set $\cup_{k'\in \mathbb{Z}^2}\partial Z'_{k',1}$). Moreover, for every $\varepsilon>0$, we have
$$\kappa\left({x'\over \varepsilon^\ell}\right)=k'\Longleftrightarrow x'\in Z'_{k',\varepsilon^\ell}.$$

\item For $k^{\prime}\in \mathcal{T}_{\varepsilon}$, the restriction of $(\widehat{\bf u}_{\varepsilon},  \widehat T_{\varepsilon}, \widehat{p}_{\varepsilon}^1)$  to $Z^{\prime}_{k^{\prime},{\varepsilon^\ell}}\times Z$ does not depend on $x^{\prime}$, whereas as a function of $z$ it is obtained from $(\widetilde{\bf u}_{\varepsilon},  \widetilde T_{\varepsilon}, \widetilde{p}_{\varepsilon}^1)$ by using the change of variables
\begin{equation}\label{CVunfolding}
\displaystyle z^{\prime}=\frac{x^{\prime}- {\varepsilon^\ell}k^{\prime}}{{\varepsilon^\ell}},\end{equation}
which transforms $Z_{k^{\prime}, {\varepsilon^\ell}}$ into $Z$. Analogously, the restriction of $\widehat p_\ep^0$ to $Z'_{k',\ep^\ell}\times Z'$ does not depend on $x'$, while as function of $z'$ it is obtained from $p_\ep^0$ by using the previous change of variables.
\end{itemize}
\end{remark}

%
Now, we deduce the following estimates for $(\widehat u_\varepsilon,\widehat T_\ep, \widehat p^0_\varepsilon, \widehat p^1_\ep)$:

\begin{lemma}\label{estimates_hat} Assume $q\in(1,3/2)$. Then,  $\widehat u_\varepsilon,\widehat T_\ep, \widehat p^0_\varepsilon$ and $\widehat p^1_\ep$ satisfy the following estimates:
  \begin{eqnarray}\medskip
 \|\widehat {\bf u}_\varepsilon\|_{L^2(\omega\times Z)^3}\leq C,&
 \|D_{z'}\widehat {\bf u}_\varepsilon\|_{L^2(\omega\times Z)^{3\times 2}}\leq C\varepsilon^{\ell-1},&
  \|\partial_{z_3}\widehat {\bf u}_\varepsilon\|_{L^2(\omega\times Z)^{3}}\leq  C,\label{estim_u_hat}\\
  \medskip
   \|\widehat T_\ep\|_{L^q(\omega\times Z)}\leq C,&
 \|\nabla_{z'}\widehat T_\ep\|_{L^q(\omega\times Z)^{2}}\leq C\varepsilon^{\ell-1},&
  \|\partial_{z_3}\widehat T_\ep\|_{L^q(\omega\times Z)}\leq  C,\label{estim_T_hat}\\
  \medskip
 \|\widehat p^0_\varepsilon\|_{L^2(\omega\times Z')}\leq C\ep^{-2},& \|\nabla_{z'}\widehat p^0_\varepsilon\|_{L^2(\omega\times Z')}\leq C\ep^{\ell-2},&  \|\widehat p^1_\varepsilon\|_{L^2(\omega\times Z')}\leq C\ep^{-1}.\label{estim_P_hat}
    \end{eqnarray}
 \end{lemma}
\begin{proof}
The proof is based in the relations between the  general relations of the norms of the unfolding unknowns and norms of the original unknowns given in  \cite[Lemma 4.9]{Anguiano_SG}:
\begin{itemize}
\item  For every $\widetilde\varphi_\varepsilon \in L^r(\widetilde\Omega_\varepsilon)^3$, $1\leq r<+\infty$,
$$\|\widehat \varphi_\varepsilon\|_{L^r(\omega\times Z)^3}=  \|\widetilde \varphi_\varepsilon\|_{L^r(\widetilde\Omega_\varepsilon)^3},$$
where $\widehat \varphi_\varepsilon$ is given by the change of variables (\ref{phihat})$_1$ (the same estimates is for (\ref{phihat})$_2$).  Similarly, for every $\zeta_\ep\in L^{r'}(\omega)$ and $\widetilde\psi_\ep \in L^{r'}(\Omega)$, the functions $\widehat \zeta_\ep$ and $\widehat \psi_\varepsilon$, given by the change of variables (\ref{phihat})$_{3,4}$ respectively satisfy
$$\|\widehat \zeta_\varepsilon\|_{L^{r'}(\omega\times Z')} = \|\zeta_\varepsilon\|_{L^{r'}(\omega)},\quad \|\widehat \psi_\varepsilon\|_{L^{r'}(\omega\times Z)} = \|\widetilde \psi_\varepsilon\|_{L^{r'}(\widetilde \Omega_\ep)}.$$
\item  For every $\widetilde \varphi_\ep\in W^{1,r}(\widetilde\Omega_\varepsilon)^3$, $1\leq r<+\infty$, the function $\widehat \varphi_\varepsilon$ given by (\ref{phihat})$_1$ belongs to $L^r(\omega;W^{1,r}(Z)^3)$, and
$$\|D_{z'} \widehat \varphi_\varepsilon\|_{L^r(\omega\times Z)^{3\times 2}} = \varepsilon^\ell  \|D_{x'}\widetilde \varphi_\varepsilon\|_{L^r(\widetilde\Omega_\varepsilon)^{3\times 2}},\quad \|\partial_{z_3} \widehat \varphi_\varepsilon\|_{L^r(\omega\times Z)^{3 }} =  \|\partial_{z_3}\widetilde \varphi_\varepsilon\|_{L^r(\widetilde\Omega_\varepsilon)^{3}},$$
$$ \|\mathbb{D}_{z'}[\widehat \varphi_\varepsilon]\|_{L^r(\omega\times Z)^{3\times 2}} = \varepsilon^\ell  \|\mathbb{D}_{x'}[\widetilde \varphi_\varepsilon]\|_{L^r(\widetilde\Omega_\varepsilon)^{3\times 2}},\quad \|\partial_{z_3}[\widehat \varphi_\varepsilon]\|_{L^r(\omega\times Z)^{3 }} = \|\partial_{z_3}[\widetilde \varphi_\varepsilon]\|_{L^r(\widetilde\Omega_\varepsilon)^{3}}.$$
\item For every $\zeta_\ep\in W^{1,r'}(\omega_\ep)^3$, $1\leq r'<+\infty$, the function $\widehat \zeta_\ep$ given by (\ref{phihat})$_{2}$ belongs to $L^{r'}(\omega;W^{1,r'}(Z')^3)$, and
$$\|\nabla_{z'}\widehat \zeta_\ep\|_{L^{r'}(\omega\times Z')^2}=\ep^{\ell}\|\nabla_{x'}\zeta_\ep\|_{L^{r'}(\omega)^2}.$$
\end{itemize}
Finally, using previous relations together with estimates for velocity, temperature and pressures given in Propositions \ref{lemma_estimates} and \ref{prop_decomposition}, we deduce the desired estimates (\ref{estim_u_hat})-(\ref{estim_P_hat}).

\end{proof}

\section{Convergence results}\label{sec:conv}
In this section, we analyze the asymptotic behavior of the functions $(\widetilde {\bf u}_\ep, \widetilde T_\ep, p_\ep^0, \widetilde p_\ep^1)$ and the corresponding unfolded functions $(\widehat {\bf u}_\ep, \widehat T_\ep, \widehat p^0_\ep, \widehat p^1_\ep)$, when $\ep$ tends to zero. It should be noted that the asymptotic behavior of the velocities and temperatures depends on the considered roughness regime, i.e.~on the parameter $\ell$. Thus, in Section 5 we will derive different limit systems, depending on the value of $\ell$.
\\

We define the following sets for $1<r<+\infty$. Let $C^\infty_{\#}(Z)$ be the space of infinitely differentiable functions in $\mathbb{R}^3$ that are $Z'$-periodic. By $L^r_{\#}(Z)$,  we denote its completion in the norm $L^r(Z)$ and by $L^r_{0,\#}(Z)$  the space of functions in $L^r_{\#}(Z)$ with zero mean value. Moreover, we introduce
\begin{equation}\label{Vz3}
\begin{array}{c}
\displaystyle
W^{1,r}_\#(Z)=\left\{\varphi\in W^{1,r}_{\rm loc}(Z)\cap L^r_\#(Z)\,:\,\int_{Z}|\nabla_z\varphi|^r\,dz<+\infty\right\},
\\
\noame
\displaystyle
V_{z_3}^r(\Theta)=\{\varphi\in L^r(\Theta)\ :\ \partial_{z_3}\varphi\in L^r(\Theta)\},\quad V_{z_3,\#}^r(\Theta)=\{\varphi\in L^r_\#(\Theta)\ :\ \partial_{z_3}\varphi\in L^q(\Theta)\}\\
\noame
\displaystyle V_{z_3,\#}^r(\omega\times Z)=\{\varphi\in L^r(\omega;L^r_\#(Z))\ :\ \partial_{z_3}\varphi\in L^r(\omega\times Z)\}.
\end{array}
\end{equation}
\begin{lemma}[Convergence of the velocity] \label{lem_conv_vel} The velocities $\widetilde {\bf u}_\ep$ (extended to $\Omega$) and $\widehat {\bf u}_\ep$ satisfy the following convergence results depending on the value of $\ell$:
\begin{itemize}
\item For $0<\ell<1$, then there exist
\begin{itemize}
\item  $\widetilde {\bf u}\in V_{z_3}^2(\Omega)^3$ where $\widetilde {\bf u}=0$ on $\Gamma_0\cup \Gamma_1$ and $\widetilde u_3\equiv 0$, such that, up to a subsequence,
\begin{equation}\label{conv_vel_tilde2}
 \widetilde {\bf u}_\ep\rightharpoonup \widetilde {\bf u}=(\widetilde {\bf u}',0)\quad\hbox{in }V_{z_3}^2(\Omega)^3,
\end{equation}
\item  $\widehat {\bf u}\in  V_{z_3,\#}^2(\omega\times Z)^3$, where $\widehat {\bf u}=0$ on $\omega\times (\widehat \Gamma_0\cup \widehat  \Gamma_1)$ and $\widehat u_3\equiv 0$, such that, up to a subsequence,
\begin{equation}\label{conv_vel_hat2}
 \widehat {\bf u}_\ep\rightharpoonup \widehat {\bf u}=(\widehat {\bf u}',0)\quad\hbox{in }V_{z_3}^2(\omega\times Z)^3,
\end{equation}
\begin{equation}\label{div_conv_z}
{\rm div}_{z'}(\widehat {\bf u}')=0\quad\hbox{in }\omega\times Z.
\end{equation}
\end{itemize}

\item For $\ell=1$, then there exist
\begin{itemize}
\item  $\widetilde {\bf u}\in V_{z_3}^2(\Omega)$ where $\widetilde {\bf u}=0$ on $\Gamma_0\cup \Gamma_1$ and $\widetilde u_3\equiv 0$, such that, up to a subsequence,
\begin{equation} \label{convStokestilde}
 \widetilde {\bf u}_\ep\rightharpoonup (\widetilde {\bf u}',0)\quad\hbox{in }V_{z_3}^2(\Omega)^3,
\end{equation}
\item  $\widehat {\bf u}\in  L^2(\omega;H^1_\#(Z))^3$, where $\widehat {\bf u}=0$ on $\omega\times (\widehat \Gamma_0\cup \widehat  \Gamma_1)$ such that, up to a subsequence,
\begin{equation} \label{convhatStokes}
 \widehat {\bf u}_\ep\rightharpoonup \widehat {\bf u}\quad\hbox{in }L^2(\omega;H^1(Z))^3,
\end{equation}
\begin{equation}\label{div_conv_z_Stokes}
{\rm div}_{z}(\widehat {\bf u})=0\quad\hbox{in }\omega\times Z.
\end{equation}
\end{itemize}
\end{itemize}
Moreover,  in every case, it holds that
\begin{equation}\label{div_conv}
{\rm div}_{x'}\left(\int_0^{h_{\rm max}}\widetilde {\bf u}'(x',z_3)\,dz_3\right)=0\quad \hbox{in }\omega,\quad\left(\int_0^{h_{\rm max}}\widetilde {\bf u}'(x',z_3)\,dz_3\right)\cdot n=0\quad \hbox{on }\partial\omega,
\end{equation}
and the following relation holds
\begin{equation}\label{relation}
\widetilde {\bf u}(x',z_3)=\int_{Z'}\widehat {\bf u}(x',z)\,dz',\quad \hbox{with}\quad \int_{Z'}\widehat u_3(x',z)\,dz'=0,
\end{equation}
and so,
\begin{equation}\label{relation2}
\int_0^{h_{\rm max}}\widetilde {\bf u}(x',z_3)\,dz_3=\int_{Z}\widehat {\bf u}(x',z)\,dz,\quad \hbox{with}\quad \int_{Z}\widehat u_3(x',z)\,dz=0.
\end{equation}
\end{lemma}
 \begin{proof}
We divide the proof in two steps, depending on the value of $\ell$.

{\it Step 1. Case $0<\ell<1$}.  Convergence (\ref{conv_vel_tilde2}) follows from  estimates for extended $\widetilde {\bf u}_\ep$ (\ref{estim_sol_dil1}) in $\Omega$, i.e.
$$\|\widetilde{\bf u}_\ep\|_{L^2(\Omega)^3}\leq C,\quad \|D_{x'}\widetilde{\bf u}_\ep\|_{L^2(\Omega)^{3\times 2}}\leq C\ep^{-1},\quad \|\partial_{z_3}\widetilde{\bf u}_\ep\|_{L^2(\Omega)^3}\leq C.$$
The proof of  $\widetilde {\bf u}=0$ on $\Gamma_0\cup \Gamma_1$, $\widetilde u_3\equiv 0$ and divergence property (\ref{div_conv}) are given in \cite{Pazanin_SG_Darcy}, so we omit it.

On the other hand, from estimates (\ref{estim_u_hat}), we have
$$\|\widehat {\bf u}_\varepsilon\|_{L^2(\omega\times Z)^3}\leq C,\quad
 \|D_{z'}\widehat {\bf u}_\varepsilon\|_{L^2(\omega\times Z)^{3\times 2}}\leq C\varepsilon^{\ell-1},\quad
  \|\partial_{z_3}\widehat {\bf u}_\varepsilon\|_{L^2(\omega\times Z)^{3}}\leq  C.$$
Since $\ep\ll \ep^\ell$, by proceeding as for $\widetilde{\bf  u}_\ep$, we can deduce convergences
$$\widehat {\bf u}_\ep\rightharpoonup \widehat {\bf u}\quad\hbox{in }V_{z_3}^2(\omega\times Z)^3,$$
and $\ep^{1-\ell}\partial_{z_i}\widehat {\bf u}_\ep$ tends to zero, $i=1,2$. Also, it holds that $\widehat {\bf u}=0$ on $\omega\times (\widehat \Gamma_0\cup \widehat \Gamma_1)$. By applying the unfolding change of variables to the divergence equation ${\rm div}_\ep(\widetilde {\bf u}_\ep)=0$ in $\widetilde\Omega^\ep$, we have that (after multiplying by $\ep$)
\begin{equation}\label{div_hat}\ep^{1-\ell}{\rm div}_{z'}(\widehat{\bf u}'_\ep)+\partial_{z_3}\widehat u_3=0\quad\hbox{in }\omega\times Z,
\end{equation}
and passing to the limit (since $\ep^{1-\ell}\to 0$ and ${\rm div}_{z'}(\widehat{\bf u}'_\ep)$ is bounded), we deduce $\partial_{z_3}\widehat u_3=0$, and so $\widehat u_3\equiv 0$. It would remain to prove the $Z'$-periodicity of $\widehat {\bf u}$ in $z'$. This can be obtain by proceeding as in
\cite[Lemma 5.4]{grau1}. This finishes the proof of convergences (\ref{conv_vel_hat2}).

Now, taking a test function $\ep^{\ell -1}\widehat \varphi \in C^1_c(\omega\times Z')$ in (\ref{div_hat}), after integrating by parts, we deduce
$$\int_{\omega\times Z}\widehat {\bf u}'\cdot \nabla_{z'}\widehat \varphi\,dx'dz=0,$$
and passing to the limit as $\ep$ tends to zero, after integrating by parts, we deduce (\ref{div_conv_z}).

The proof of relation (\ref{relation}) is similar to the step 3 of the proof of \cite[Lemma 5.4]{Anguiano_SG}, so we omit it.
\\

 {\it Step 2. Case $\ell=1$}.   The proof of convergence (\ref{convStokestilde}) is similar to the previous case, so we omit it. On the other hand, from estimates (\ref{estim_u_hat})  with $\ell=1$, we have
$$\|\widehat {\bf u}_\varepsilon\|_{L^q(\omega\times Z)^3}\leq C,\quad
 \|D_{z'}\widehat {\bf u}_\varepsilon\|_{L^q(\omega\times Z)^{3\times 2}}\leq C,\quad
  \|\partial_{z_3}\widehat {\bf u}_\varepsilon\|_{L^q(\omega\times Z)^{3}}\leq  C.$$
Then, we deduce convergences
$$\widehat {\bf u}_\ep\rightharpoonup \widehat {\bf u}\quad\hbox{in }L^2(\omega;H^1(Z))^3,$$
and also, it holds that $\widehat {\bf u}=0$ on $\omega\times (\widehat \Gamma_0\cup \widehat \Gamma_1)$. By applying the unfolding change of variables to the divergence equation and multiplying by $\ep$, we have that
\begin{equation}\label{div_hatStokes} {\rm div}_{z}(\widehat{\bf u}_\ep)=0\quad\hbox{in }\omega\times Z,
\end{equation}
and passing to the limit, we deduce divergence condition (\ref{div_conv_z_Stokes}). The $Z'$-periodicity of $\widehat {\bf u}$ in $z'$  and the rest of properties can be deduced as in the previous case, so we omit it.

 \end{proof}

\begin{lemma}[Convergences of temperatures] \label{lem_conv_temp}Assume $q\in(1,3/2)$. The  temperatures $\widetilde T_\ep$ (extended to $\Omega$) and $\widehat T_\ep$ satisfy the following convergence results depending on the value of $\ell$:
\begin{itemize}
\item For $0<\ell<1$, then there exist
\begin{itemize}
\item  $\widetilde T\in V_{z_3}^q(\Omega)$ where $\widetilde T=0$ on $\Gamma_0\cup \Gamma_1$, such that, up to a subsequence,
\begin{equation}\label{conv_vel_tilde2temp}
 \widetilde T_\ep\rightharpoonup \widetilde T\quad\hbox{in }V_{z_3}^q(\Omega),
\end{equation}
\item  $\widehat T\in  V_{z_3,\#}^q(\omega\times Z)$, where $\widehat T=0$ on $\omega\times (\widehat \Gamma_0\cup \widehat  \Gamma_1)$, such that, up to a subsequence,
\begin{equation}\label{conv_vel_hat2temp}
\widehat T_\ep\rightharpoonup \widehat T\quad\hbox{in }V_{z_3}^q(\omega\times Z),
\end{equation}
\end{itemize}

\item For $\ell=1$, then there exist
\begin{itemize}
\item  $\widetilde T\in V_{z_3}^q(\Omega)$ where $\widetilde T=0$ on $\Gamma_0\cup \Gamma_1$, such that, up to a subsequence,
\begin{equation} \label{convStokestildetemp}
 \widetilde T_\ep\rightharpoonup \widetilde T\quad\hbox{in }V_{z_3}^q(\Omega),
\end{equation}
\item  $\widehat T\in  L^q(\omega;W^{1,q}_\#(Z))$, where $\widehat T=0$ on $\omega\times (\widehat \Gamma_0\cup \widehat  \Gamma_1)$ such that, up to a subsequence,
\begin{equation} \label{convhatStokestemp}
 \widehat T_\ep\rightharpoonup \widehat T\quad\hbox{in }L^q(\omega;W^{1,q}(Z)),
\end{equation}
\end{itemize}
\end{itemize}
Moreover,  in every case,  the following relation holds
\begin{equation}\label{relationtemp}
\widetilde T(x',z_3)=\int_{Z'}\widehat T(x',z)\,dz',
\end{equation}
and so,
\begin{equation}\label{relation2temp}
\int_0^{h_{\rm max}}\widetilde T(x',z_3)\,dz_3=\int_{Z}\widehat T(x',z)\,dz.
\end{equation}
\end{lemma}
\begin{proof}
The proof is similar to the proof of velocity convergences, just changing the Sobolev spaces with index $q$ instead of $2$, so we omit it.

\end{proof}
\begin{lemma}[Convergences of pressures]\label{lem_conv_press} There exist $\widetilde p\in L^{2}_0(\omega)\cap H^1(\omega)$, independent of $z$ and $\widehat p_1\in L^{2}(\omega;L^{2}_\#(Z))$ such that
\begin{equation}\label{pressuresH}
\ep^2 p_\ep^0\rightharpoonup  p\quad\hbox{in }H^1(\omega),\quad \ep^2 p_\ep^0\to \widetilde p\quad\hbox{in }L^{2}(\omega),\quad \ep^2\widehat p_\ep^0\to \widetilde p\quad\hbox{in } L^2(\omega\times Z),
\end{equation}
\begin{equation}\label{pressuresp1}
\ep\widehat p_1^\ep\rightharpoonup\widehat p_1\quad\hbox{in }L^{2}(\omega;L^{2}(Z)).
\end{equation}
\end{lemma}
\begin{proof}
From estimates for $p_0^\ep$ and $\widehat p_\ep^0$, and the classical compactness result for the unfolding method for a bounded sequence in $H^1$, see for instance \cite{Cioran-book}, we obtain convergences (\ref{pressuresH}). Estimates for $\widehat p_1^\ep$ given in Proposition \ref{prop_decomposition} imply the existence of $\widehat p_1\in L^2(\omega; L^2_{\#}(Z))$ such that up to a subsequence, convergence (\ref{pressuresp1}) holds.

Finally, since $\widetilde p_\ep$ has mean value zero in $\widetilde\Omega_\ep$, from the decomposition of the pressure and the unfolding change of variables, we deduce
$$0=\int_{\omega\times Z}\widehat p_\ep^0\,dx'dz'+\int_{\omega\times Z}\widehat p^1_\ep\,dx'dz.$$
Considering the convergence of $\widehat p^0_\ep$ to $\widetilde p$, that $\widehat p^1_\ep$ tends to zero, and that $\widetilde p$ does not depend on $z'$, we obtain
$$\int_\omega \widetilde p\,dx'=0,$$
and then, $\widetilde p$ has mean value zero in $\omega$.

\end{proof}
We finish this section by giving a variational formulation for the unfolded functions. To do this, we consider the following functions:
\begin{itemize}
\item According to Lemma \ref{lem_conv_vel}, we choose   ${\bf v}(x',z)\in \mathcal{D}(\omega; C^\infty_{\#}(Z)^3)$  and boundary values   ${\bf v}(x',z)=0$ on $\omega \times (\widehat \Gamma_0\cup \widehat \Gamma_1)$ and ${\rm div}_{x'}(\int_Z {\bf v}'\,dz)=0$ in $\omega$ and $(\int_{Z}{\bf v}'\,dz )\cdot n=0$ on $\partial\omega$.  Moreover, in the case $0<\ell<1$, we assume $v_3\equiv 0$ and ${\rm div}_{z'}({\bf v}')=0$ in $\omega\times Z$, while in the case $\ell=1$, we assume ${\rm div}_{z}({\bf v})=0$ in $\omega\times Z$.

\item According to Lemma \ref{lem_conv_temp}, we choose $\zeta(x',z)\in \mathcal{D}(\omega; C^\infty_{\#}(Z)^3)$  and boundary values   $\zeta(x',z)=0$ on $\omega \times (\widehat \Gamma_0\cup \widehat \Gamma_1)$.
\end{itemize}

Thus, taking $\ep^2 {\bf v}(x',x'/\varepsilon^\ell,z_3)$ as a test function in (\ref{formvar_dil})$_1$ and taking into account the decomposition of the pressure, we get
$$
\begin{array}{l}
\displaystyle 2\mu_{\rm eff}\int_{\widetilde \Omega^\ep}\!\!\mathbb{D}_{\ep}[\widetilde  {\bf u}_\ep]:(\mathbb{D}_{x'}[  {\bf v}])+\ep^{2-\ell}\mathbb{D}_{z'}[   {\bf v}]+\ep\partial_{z_3}[{\bf v}])\,dx'dz_3\!
+\!{\mu\over  K}\int_{\widetilde \Omega^\ep}\!\!\widetilde  {\bf u}_\ep\cdot  {\bf v}\,dx'z_3\\
\noame
\displaystyle -\ep^2\int_{\widetilde \Omega_\ep}(p^0_\ep+\widetilde p^1_\varepsilon) \left({\rm div}_{x'}({\bf v}')+\ep^{-\ell} {\rm div}_{z'}({\bf v}')+\ep^{-1}\partial_{z_3}v_3\right)\,dx'dz_3=\int_{\widetilde \Omega^\ep} {\bf f}'\cdot    {\bf v}'\,dx'dz_3.
\end{array}$$
Similarly, by taking $\zeta(x',x'/\varepsilon^\ell,z_3)$ as a test function in (\ref{formvar_dil})$_3$, we deduce
 $$
\begin{array}{l}
 \displaystyle \ep^2\int_{\widetilde \Omega^\ep}  k\nabla_{x'}\widetilde  T_\ep\cdot \nabla_{x'}  \zeta\,dx'dz_3+\ep^{2-\ell}\int_{\widetilde \Omega^\ep} k\nabla_{x'}\widetilde  T_\ep\cdot \nabla_{z'}  \zeta\,dx'dz_3
 +\int_{\widetilde \Omega^\ep}
k\partial_{z_3}\widetilde  T_\ep\,\partial_{z_3}\zeta\,dx'dz_3
 \\
 \noame
 \displaystyle
 ={\mu\over K} \int_{\widetilde \Omega^\ep} |\widetilde  {\bf u}_\ep|^2  \zeta\,dx'dz_3+2\mu_{\rm eff}\ep^{2}\int_{\widetilde \Omega^\ep}|\mathbb{D}_{\ep}[\widetilde  {\bf u}_\ep]|^2   \zeta\,dx'dz_3+\int_{\omega}   b\,  \zeta\,dx'.
\end{array}$$
By using the change of variables given in Remark \ref{remarkCV} in previous variational formulations, we respectively obtain
\begin{equation}\label{form_var_hat}\begin{array}{l}
\displaystyle
\displaystyle 2\mu_{\rm eff}\int_{\omega\times Z}(\ep^{-\ell}\mathbb{D}_{z'}[\widehat  {\bf u}_\ep]+\ep^{-1}\partial_{z_3}[\widehat{\bf u}_\ep]):(\ep^{2-\ell}\mathbb{D}_{z'}[ {\bf v}])+\ep\partial_{z_3}[{\bf v}])\,dx'dz\!
+\!{\mu\over  K} \int_{\omega\times Z}\!\!\widehat  {\bf u}_\ep\cdot   {\bf v}\,dx'dz\\
\noame
\displaystyle -\int_{\omega\times Z}\ep^2(\widehat p_\ep^0+\widehat  p^\varepsilon_1) \left({\rm div}_{x'}({\bf v}')+\ep^{-\ell} {\rm div}_{z'}({\bf v}')+\ep^{-1}\partial_{z_3}v_3\right)\,dx'dz=\int_{\omega\times Z} {\bf f}'\cdot   {\bf v}'\,dx'dz + O_\ep,
\end{array}
\end{equation}
and
\begin{equation}\label{form_var_hat_temp}
\begin{array}{l}
 \displaystyle \ep^{2(1-\ell)}\int_{\omega\times Z} k\nabla_{z'}\widehat  T_\ep\cdot \nabla_{z'}  \zeta\,dx'dz
 +\int_{\omega\times Z}
k\partial_{z_3}\widehat  T_\ep\,\partial_{z_3}\zeta\,dx'dz
 \\
 \noame
 \displaystyle
 ={\mu\over K} \int_{\omega\times Z} |\widehat  {\bf u}_\ep|^2  \zeta\,dx'dz+2\mu_{\rm eff}\ep^{2}\int_{\omega\times Z}|\ep^{-\ell}\mathbb{D}_{z'}[\widehat  {\bf u}_\ep]+\ep^{-1}\partial_{z_3}[\widehat {\bf u}_\ep]|^2  \zeta\,dx'dz+ \int_{\omega\times Z'}   b\,  \zeta\,dx'dz'+O_\ep,
\end{array}
\end{equation}
where $O_\varepsilon$ is a generic real sequence depending on $\varepsilon$, that can change from line to line, and devoted to tend to zero. Finally, depending on the value of $\ell$, we have the following variational formulations:
\begin{itemize}
\item In the case $0<\ell<1$, taking into account that $v_3\equiv 0$ and ${\rm div}_{z'}({\bf v}')=0$ in $\omega\times Z$, then (\ref{form_var_hat}) and (\ref{form_var_hat_temp}) read as follows
\begin{equation}\label{form_var_hat_subcritical}\begin{array}{l}
\displaystyle
\displaystyle 2\mu_{\rm eff}\int_{\omega\times Z}\ep^{2(1-\ell)}\mathbb{D}_{z'}[\widehat  {\bf u}_\ep']: \mathbb{D}_{z'}[  {\bf v}']\,dx'dz+2\mu_{\rm eff}\int_{\omega\times Z}\partial_{z_3}[\widehat{\bf u}_\ep']:\partial_{z_3}[ {\bf v}']\,dx'dz
\displaystyle+\!{\mu\over  K} \int_{\omega\times Z}\!\!\widehat  {\bf u}_\ep'\cdot   {\bf v}'\,dx'dz\\
\noame
\displaystyle -\int_{\omega\times Z}\ep^2(\widehat p_\ep^0+\widehat  p^\varepsilon_1) {\rm div}_{x'}({\bf v}')\,dx'dz=\int_{\omega\times Z}{\bf f}'\cdot   {\bf v}'\,dx'dz + O_\ep,
\\
\\
 \displaystyle \ep^{2(1-\ell)}\int_{\omega\times Z} k\nabla_{z'}\widehat  T_\ep\cdot \nabla_{z'}  \zeta\,dx'dz
 +\int_{\omega\times Z}
k\partial_{z_3}\widehat  T_\ep\,\partial_{z_3}\zeta\,dx'dz
 \\
 \noame
 \displaystyle
 ={\mu\over K} \int_{\omega\times Z} |\widehat  {\bf u}_\ep|^2  \zeta\,dx'dz+2\mu_{\rm eff}\int_{\omega\times Z}|\ep^{1-\ell}\mathbb{D}_{z'}[\widehat  {\bf u}_\ep]+ \partial_{z_3}[\widehat {\bf u}_\ep]|^2  \zeta\,dx'dz+ \int_{\omega\times Z'}   b\,  \zeta\,dx'dz'+O_\ep.
\end{array}
\end{equation}

\item In the case $\ell=1$, taking into account that ${\rm div}_z({\bf v})=0$ in $\omega\times Z$, then (\ref{form_var_hat}) and (\ref{form_var_hat_temp}) read as follows
\begin{equation}\label{form_var_hat_templ1}\begin{array}{l}
\displaystyle
\displaystyle 2\mu_{\rm eff}\int_{\omega\times Z} \mathbb{D}_{z}[\widehat  {\bf u}_\ep]:\mathbb{D}_{z}[ {\bf v}]\,dx'dz\!
+\!{\mu\over  K} \int_{\omega\times Z}\!\!\widehat  {\bf u}_\ep\cdot   {\bf v}\,dx'dz -\int_{\omega\times Z}\ep^2(\widehat p_\ep^0+\widehat  p^\varepsilon_1)\,{\rm div}_{x'}({\bf v}')\,dx'dz\\
\noame
\displaystyle=\int_{\omega\times Z} {\bf f}'\cdot {\bf v}'\,dx'dz + O_\ep,
\\
\\
 \displaystyle  \int_{\omega\times Z} k\nabla_{z}\widehat  T_\ep\cdot \nabla_{z}  \zeta\,dx'dz
  ={\mu\over K} \int_{\omega\times Z} |\widehat  {\bf u}_\ep|^2  \zeta\,dx'dz+2\mu_{\rm eff} \int_{\omega\times Z}| \mathbb{D}_{z}[\widehat  {\bf u}_\ep]|^2  \zeta\,dx'dz+ \int_{\omega\times Z'}   b\,  \zeta\,dx'dz'+O_\ep.
\end{array}
\end{equation}
\end{itemize}

\section{Limit problems}
\subsection{The case $0<\ell<1$}
We start by proving that the limit functions $(\widehat{\bf u}, \widehat T, \widetilde p)$ satisfies a two-pressure problem. Then, we derive the main result for the case $0<\ell<1$ by providing the expressions for velocity $(\widetilde {\bf u}, \widetilde T)$ and the Reynolds problem satisfied by $\widetilde p$.
\begin{theorem}[Limit model]\label{thm_limit_0ell1} The limit functions $\widehat{\bf  u}$, $\widehat T$ and $\widetilde  p$ given in Lemmas \ref{lem_conv_vel}, \ref{lem_conv_temp} and \ref{lem_conv_press} satisfy
\begin{equation}\label{limit_model_0ell1}
\left\{\begin{array}{rl}
\displaystyle
-\mu_{\rm eff} \partial_{z_3}^2 \widehat {\bf  u}' +{\mu\over K}\widehat {\bf u}'+\nabla_{z'}\widehat \pi={\bf f}'(x') - \nabla_{x'}\widetilde  p(x')&\hbox{ in }\omega\times Z,\\
\noame
\displaystyle {\rm div}_{z'}\left(\int_0^{h(z')}\widehat {\bf u}'\,dz_3\right)=0&\hbox{ in }\omega\times Z',\\
\noame
\displaystyle-k\partial^2_{z_3}\widehat T={\mu\over K}|\widehat {\bf  u}'|^2+\mu_{\rm eff}|\partial_{z_3}\widehat {\bf  u}'|^2 &\hbox{ in }\omega\times Z,\\
\noame

\widehat {\bf  u}'=0&\hbox{ on } \omega\times (\widehat \Gamma_0\cup \widehat \Gamma_1),\\
\noame
\widehat T=0&\hbox{ on }    \omega\times \widehat \Gamma_1,\\
\noame
-k\partial_{z_3}{ \widehat T}=b&\hbox{ on }\omega\times \widehat \Gamma_0,\\
\noame
\widehat \pi\in L^{2}(\omega;L^{2}_{0,\#}(Z')).
\end{array}\right.
\end{equation}
\end{theorem}
\begin{proof}
 To prove system (\ref{limit_model_0ell1}), we observe that conditions (\ref{limit_model_0ell1})$_{2, 4, 5}$ are consequence of Lemmas \ref{lem_conv_vel} and \ref{lem_conv_temp}. We divide the proof in four steps. First, we pass to the limit in the variational formulation for velocity (\ref{form_var_hat_subcritical})$_1$. Then, in second and third steps, we prove strong convergences of velocity, which is used to pass to the limit  in  the variational formulation for temperature (\ref{form_var_hat_subcritical})$_2$ in the fourth step.

{\it Step 1.} Let us pass to the limit in the variational formulation for velocity (\ref{form_var_hat_subcritical})$_1$:
\begin{itemize}
\item First, second and third terms.  According to estimate (\ref{estim_u_hat})$_2$ and $\ep^{1-\ell}\to 0$, we deduce
that
$$\left|2\mu_{\rm eff}\int_{\omega\times Z}\ep^{2(1-\ell)}\mathbb{D}_{z'}[\widehat  {\bf u}_\ep']: \mathbb{D}_{z'}[  {\bf v}']\,dx'dz\right|\leq C\ep^{1-\ell}\to 0,$$
and so, from convergence of velocity given in Lemma \ref{lem_conv_vel}, we get
$$\begin{array}{l}
\displaystyle 2\mu_{\rm eff}\int_{\omega\times Z}\ep^{2(1-\ell)}\mathbb{D}_{z'}[\widehat  {\bf u}_\ep']: \mathbb{D}_{z'}[  {\bf v}']\,dx'dz+2\mu_{\rm eff}\int_{\omega\times Z}\partial_{z_3}[\widehat{\bf u}_\ep']:\partial_{z_3}[ {\bf v}']\,dx'dz
\displaystyle+\!{\mu\over  K} \int_{\omega\times Z}\!\!\widehat  {\bf u}_\ep'\cdot   {\bf v}'\,dx'dz\\
\noame
\displaystyle
\to 2\mu_{\rm eff}\int_{\omega\times Z}\partial_{z_3}[\widehat{\bf u}']:\partial_{z_3}[ {\bf v}']\,dx'dz
\displaystyle+\!{\mu\over  K} \int_{\omega\times Z}\!\!\widehat  {\bf u}'\cdot   {\bf v}'\,dx'dz.
\end{array}$$
\item Fourth term. According to estimate (\ref{estim_p0p1})$_3$, we deduce
$$\left|\int_{\omega\times Z}\ep^2\widehat  p^\varepsilon_1 {\rm div}_{x'}({\bf v}')\,dx'dz\right|\leq C\ep\to 0,$$
and so, from convergence of pressure given in Lemma \ref{lem_conv_press}, we deduce
$$\int_{\omega\times Z}\ep^2(\widehat p_\ep^0+\widehat  p^\varepsilon_1) {\rm div}_{x'}({\bf v}')\,dx'dz\to \int_{\omega\times Z}\widetilde p\,{\rm div}_{x'}({\bf v}')\,dx'dz.$$
Observe that from condition ${\rm div}_{x'}(\int_Z {\bf v}'\,dz)=0$ in $\omega$ and since $\widetilde p$ does not depend on $z$, we get
$$\int_{\omega\times Z}\widetilde p\,{\rm div}_{x'}({\bf v}')\,dx'dz=\int_{\omega}\widetilde p\,{\rm div}_{x'}\left(\int_Z {\bf v}'dz\right)\,dx'=0.$$
\end{itemize}
Thus, from previous convergences, and taking into account
\begin{equation}\label{uvpartial}\partial_{z_3}[\widehat {\bf u}']: \partial_{z_3} [{\bf v}']={1\over 2}\partial_{z_3}\widehat {\bf u}' \cdot \partial_{z_3} {\bf v}',
\end{equation}
the limit of the variational formulation (\ref{form_var_hat_subcritical})$_1$ is given by
 \begin{equation}\label{var_limit_problem_0ell1}
 \mu_{\rm eff}\int_{\omega\times Z}\partial_{z_3}\widehat{\bf u}':\partial_{z_3} {\bf v}'\,dx'dz
\displaystyle+\!{\mu\over  K} \int_{\omega\times Z}\!\!\widehat  {\bf u}'\cdot   {\bf v}'\,dx'dz=\int_{\omega\times Z}{\bf f}'\cdot   {\bf v}'\,dx'dz,
 \end{equation}
 which, by density, holds for every ${\bf v}'\in \mathcal{V}$ with
$$\mathcal{V}=\left\{\begin{array}{l}
\displaystyle {\bf v}'\in V_{z_3,\#}^2(\omega\times Z')\ :\ {\rm div}_{z'}({\bf v}')=0\quad\hbox{in }\omega\times Z\\
\noame
\displaystyle
{\rm div}_{x'}\left(\int_Z{\bf v}'\,dz\right)=0\quad\hbox{in }\omega,\quad \left(\int_Z{\bf v}'\,dz\right)\cdot n=0\quad\hbox{on }\partial\omega
\end{array}
\right\}.
$$
Finally, reasoning as in \cite[Lemma 1.5]{Allaire}, the orthogonal of $\mathcal{V}$ is made of gradients of the form $\nabla_{x'}\widetilde\pi(x')+\nabla_{z'}\widehat \pi(x',z')$, with $\widetilde \pi\in L^{2}_0(\omega)$ and $\widehat\pi(x',z')\in L^{2}(\omega;L^{2}_{0,\#}(Z'))$. Thus, integrating by parts, the limit variational formulation (\ref{var_limit_problem_0ell1}) is equivalent to the two-pressures reduced Stokes system (\ref{limit_model_0ell1})$_1$. It remains to prove that $\widetilde \pi$ coincides with pressure $\widetilde p$. This can be easily done by passing to the limit similarly as above and considering the test function ${\bf v}$, which is divergence-free only in $z'$, and by identifying limits.
\\

{\it Step 2.} We prove the following property
\begin{equation}\label{prevprop_sub}\begin{array}{l}
 \displaystyle
 \lim_{\ep\to 0}\left({\mu\over K}\int_{\omega\times Z}|\widehat {\bf u}_\ep|^2\,dx'dz+2\mu_{\rm eff}\int_{\omega \times Z}| \ep^{1-\ell}\mathbb{D}_{z'}[\widehat {\bf u}_\ep]+ \partial_{z_3}[\widehat {\bf u}_\ep]|^2 \,dx'dz\right)\\
 \noame
 \displaystyle
 ={\mu\over K}\int_{\omega\times Z}|\widehat {\bf u}'|^2 \,dx'dz+2\mu_{\rm eff}\int_{\omega\times Z}|\partial_{z_3}[\widehat {\bf u}']|^2 \,dx'dz\,.
 \end{array}
 \end{equation}
We take $\widehat {\bf u}_\ep$ as test function in  (\ref{form_var_hat_subcritical})$_1$. From the fact that $\ep^{-\ell}{\rm div}_{z'}(\widehat {\bf u}_\ep')+ \ep^{-1}\partial_{z_3}\widehat u_{3,\ep}=0$, using strong convergence of $\ep^2\widehat p_\ep^0$ and convergence of velocity,  we get
\begin{equation}\label{prevpro111sub}
\begin{array}{l}
 \displaystyle
\lim_{\ep\to 0}\left(2\mu_{\rm eff}\int_{\omega \times Z} |\ep^{1-\ell}\mathbb{D}_{z'}[\widehat {\bf u}_\ep]|^2\,dx'dz+2\mu_{\rm eff}\int_{\omega \times Z} |\partial_{z_3}[\widehat {\bf u}_\ep]|^2\,dx'dz
+{\mu\over K}\int_{\omega\times Z}|\widehat {\bf u}^\ep|^2\,dx'dz\right)\\
\noame
\displaystyle=  \displaystyle\int_{\omega\times Z}\widetilde p\, {\rm div}_{x'}(\widehat {\bf u}')\,dx'dz+\int_{\omega\times Z}{\bf f}'\cdot \widehat {\bf u}'\,dx'dz=\int_{\omega\times Z}{\bf f}'\cdot \widehat {\bf u}'\,dx'dz,
\end{array}
\end{equation}
because $\int_{\omega\times Z}\widetilde p\, {\rm div}_{x'}(\widehat {\bf u}')\,dx'dz=0$, because $\widetilde p$ does not depend on $z$ and ${\rm div}_{x'}(\int_Z\widehat {\bf u}'dz)=0$.

Now, we take $\widehat {\bf u}'$ as test function in (\ref{var_limit_problem_0ell1}), and we deduce
\begin{equation}\label{limituhat2sub}\begin{array}{rl}
\displaystyle
\int_{\omega\times Z}{\bf f}'\cdot \hat{\bf u}'\,dx'dz=&\displaystyle 2\mu_{\rm eff}\int_{\omega\times Z} |\partial_{z_3}[\widehat {\bf u}']|^2\,dx'dz
+ \displaystyle {\mu\over  K}\int_{\omega\times\times Z}  |\widehat {\bf u}'|^2\,dx'dz,
 \end{array}
 \end{equation}
 where property (\ref{uvpartial}) is used.
Then, from (\ref{prevpro111sub}) and (\ref{limituhat2sub}), we get (\ref{prevprop_sub}).
\\

 {\it Step 3.} We prove the strong convergences
\begin{equation}\label{strongconvsub}\lim_{\ep\to 0}\int_{\omega\times Z} |\widehat {\bf u}_\ep-(\widehat {\bf u}',0)|^2\,dx'dz=0,\quad \lim_{\ep\to 0}\int_{\omega\times Z} \left|\ep^{1-\ell}\mathbb{D}_{z'}[\widehat {\bf u}_\ep]+\partial_{z_3}[\widehat {\bf u}_\ep]-\partial_{z_3}[\widehat {\bf u}']\right|^2\,dx'dz=0.
\end{equation}
To prove (\ref{strongconvsub}), it is enough to prove that
 $$E_\ep:= 2\mu_{\rm eff}\int_{\omega\times Z} \left| \ep^{1-\ell}\mathbb{D}_{z'}[\widehat {\bf u}_\ep]+\partial_{z_3}[\widehat {\bf u}_\ep]-\partial_{z_3}[\widehat {\bf u}']\right|^2\,dx'dz+{\mu\over K}\int_{\omega\times Z} |\widehat {\bf u}_\ep-(\widehat {\bf u}',0)|^2\,dx'dz\to 0.$$
 Developing the expression of $E_\ep$, we have
 $$\begin{array}{rl}
 E_\ep =&\displaystyle
 2\mu_{\rm eff}\int_{\omega\times Z} \left|\ep^{1-\ell}\mathbb{D}_{z'}[\widehat {\bf u}_\ep]+\partial_{z_3}[\widehat {\bf u}_\ep]\right|^2\,dx'dz+{\mu\over K}\int_{\omega\times Z} |\widehat {\bf u}_\ep|^2\,dx'dz\\
 \noame
 &\displaystyle
 +2\mu_{\rm eff}\int_{\omega\times Z} \left|\partial_{z_3}[\widehat {\bf u}']\right|^2\,dx'dz-4\mu_{\rm eff}\int_{\omega\times Z} \left(\ep^{1-\ell}\mathbb{D}_{z'}[\widehat {\bf u}_\ep]+\partial_{z_3}[\widehat {\bf u}_\ep]\right):\partial_{z_3}[\widehat {\bf u}'] \,dx'dz\\
 \noame
 &\displaystyle
 +{\mu\over K}\int_{\omega\times Z} |\widehat {\bf u}'|^2\,dx'dz-2{\mu\over K}\int_{\omega\times Z} \widehat {\bf u}_\ep'\cdot \widehat {\bf u}'\,dx'dz.
 \end{array}$$
 By using property (\ref{prevprop_sub}), convergence $\ep^{1-\ell}\to 0$   and convergences of velocity, we get
 $$\begin{array}{rl}
 E_\ep \to &\displaystyle
 2\mu_{\rm eff}\int_{\omega\times Z} \left|\partial_{z_3}[\widehat {\bf u}']\right|^2\,dx'dz+{\mu\over K}\int_{\omega\times Z} |\widehat {\bf u}'|^2\,dx'dz\\
 \noame
 &\displaystyle
 +2\mu_{\rm eff}\int_{\omega\times Z} \left|\partial_{z_3}[\widehat {\bf u}']\right|^2\,dx'dz-4\mu_{\rm eff}\int_{\omega\times Z} |\partial_{z_3}[\widehat {\bf u}']|^2\,dx'dz\\
 \noame
 &\displaystyle
 +{\mu\over K}\int_{\omega\times Z} |\widehat {\bf u}'|^2\,dx'dz-2{\mu\over K}\int_{\omega\times Z} |\widehat {\bf u}'|^2\,dx'dz=0,
 \end{array}$$
 which proves the strong convergences (\ref{strongconvsub}). \\

{\it Step 4.} Let us pass to the limit in  the temperature variational formulation (\ref{form_var_hat_subcritical})$_2$:
\begin{itemize}
\item First and second terms. By using estimate (\ref{estim_T_hat}) and $\ep^{1-\ell}\to 0$, we deduce
$$\left|\ep^{2(1-\ell)}\int_{\omega\times Z} k\nabla_{z'}\widehat  T_\ep\cdot \nabla_{z'}  \zeta\,dx'dz
\right|\leq C\ep^{1-\ell}\to 0,$$
and so, from convergence of temperature given in Lemma  \ref{lem_conv_temp}, we get
 $$\ep^{2(1-\ell)}\int_{\omega\times Z} k\nabla_{z'}\widehat  T_\ep\cdot \nabla_{z'}  \zeta\,dx'dz
 +\int_{\omega\times Z}
k\partial_{z_3}\widehat  T_\ep\,\partial_{z_3}\zeta\,dx'dz\to \int_{\omega\times Z}
k\partial_{z_3}\widehat T\,\partial_{z_3}\zeta\,dx'dz.$$

\item Third and fourth terms. From the strong convergences (\ref{strongconvsub}), we deduce
$$\begin{array}{l}\displaystyle
{\mu\over K} \int_{\omega\times Z} |\widehat  {\bf u}_\ep|^2  \zeta\,dx'dz+2\mu_{\rm eff} \int_{\omega\times Z}|\ep^{1-\ell}\mathbb{D}_{z'}[\widehat  {\bf u}_\ep]+ \partial_{z_3}[\widehat {\bf u}_\ep]|^2  \zeta\,dx'dz\\
\noame
\displaystyle
\to 2\mu_{\rm eff}\int_{\omega\times Z}|\partial_{z_3}[\widehat{\bf u}']|^2\,dx'dz
\displaystyle+\!{\mu\over  K} \int_{\omega\times Z}|\widehat  {\bf u}'|^2\,dx'dz.
\end{array}$$
\end{itemize}
Thus, from previous convergences and relation (\ref{uvpartial}),  the limit of the variational formulation (\ref{form_var_hat_subcritical})$_2$ is given by
 \begin{equation}\label{var_limit_problem_temp_0ell1}\int_{\omega\times Z}
k\partial_{z_3}\widehat T\,\partial_{z_3}\zeta\,dx'dz=\mu_{\rm eff}\int_{\omega\times Z}|\partial_{z_3}\widehat{\bf u}'|^2\,dx'dz
\displaystyle+\!{\mu\over  K} \int_{\omega\times Z}|\widehat  {\bf u}'|^2\,dx'dz+\int_{\omega\times Z'}   b\,  \zeta\,dx'dz',
\end{equation}
which, by density, holds for every $\zeta\in V_{z_3,\#}^q(\omega\times Z)$. Thus, integrating by parts, the limit variational formulation (\ref{var_limit_problem_temp_0ell1}) is equivalent to the two-pressures reduced Stokes system (\ref{limit_model_0ell1})$_{3,6}$.

\end{proof}
\begin{corollary} Consider $(\widehat u, \widetilde p, \widehat T)$, with $\widehat u_3=0$, satisfying system (\ref{limit_model_0ell1}). Let us define the average velocity and temperature by
\begin{equation}\label{averages0ell1}
{\bf V}_{\rm av}(x')=\int_0^{h_{\rm max}}\widetilde {\bf u}(x',z_3)\,dz_3,\quad T_{\rm av}(x')=\int_0^{h_{\rm max}}\widetilde T(x',z_3)\,dz_3,
\end{equation}
and a parameter $M=\sqrt{\mu\over K\mu_{\rm eff}}>0$.  We have the following:
\begin{itemize}
\item[(i)] The average velocity ${\bf V}_{\rm av}$ is given by
 \begin{equation}\label{Expressions0ell1}
{\bf V}_{\rm av}'(x',z)=-{K\over \mu}A_{M}({\bf f}'(x')-\nabla_{x'}\widetilde p(x')),\quad V_{{\rm av},3}\equiv 0,\quad \hbox{in }\omega,
\end{equation}
where $M=\sqrt{\mu\over K\mu_{\rm eff}}>0$ and $A_{M}\in \mathbb{R}^{2\times 2}$ is given by
\begin{equation}\label{AM0}
\begin{array}{l}(A_{M})_{ij}=\displaystyle -\int_{Z'}  \left({2\over M}{e^{Mh(z')}-e^{-Mh(z')}-2\over e^{Mh(z')}-e^{-Mh(z')}}- h(z')\right)(e_i+\nabla_{z'}\pi^i(z'))e_j\,dz',\quad i,j=1,2.
\end{array}
\end{equation}
Here, $\pi^i\in H^1_\#(Z')$, $i=1,2$ is the unique solution of
\begin{equation}\label{reynoldspi}\begin{array}{rl}
\displaystyle -{\rm div}_{z'}\left\{\left({2\over M}{e^{Mh(z')}-e^{-Mh(z')}-2\over e^{Mh(z')}-e^{-Mh(z')}}- h(z')\right)(e_i+\nabla_{z'}\pi^i(z'))\right\}=0&\hbox{in }Z'.
\end{array}
\end{equation}
\item[(ii)] Pressure $\widetilde p \in L^2_0(\omega)\cap H^1(\omega)$ is the unique solution of the Reynolds equation:
\begin{equation}\label{ReynoldsP}
\left\{ \begin{array}{l}\displaystyle -{\rm div}_{x'}
\left({K\over \mu}A_{M}(\nabla_{x'}\widetilde p(x')-{\bf f}'(x'))\right)=0,\quad \hbox{in }\omega, \\
\noame
\displaystyle
\left({K\over \mu}A_{M}(\nabla_{x'}\widetilde p(x')-{\bf f}'(x'))\right)\cdot n=0,\quad \hbox{on }\partial\omega.
\end{array}\right.
\end{equation}

\item[(iii)] The average temperature $T_{\rm av}$ is given by
\begin{equation}\label{ExpT_thm0ell1}
\begin{array}{l}\displaystyle  T_{\rm av}(x')= \int_{Z}\widehat T(x',z)\,dz,\\
\end{array}\end{equation}
with
\begin{equation}\label{That}
\begin{array}{rl}
\displaystyle
 \medskip
\widehat T(x',z)=&\displaystyle -{b\over k} -{\mu\over   k\,K}\left(\int_{z_3}^1 \int_{\xi}^1|\widehat{\bf u}'(x',z',\tau)|^2\,d\tau\,d\xi-z_3\int_{0}^1| \widehat{\bf u}'(x',z',\xi)|^2\,d\xi\right)\\
\noame
&\displaystyle -{\mu_{\rm eff}\over k}\left(\int_{z_3}^1\int_\xi^1|\partial_{z_3}\widehat  {\bf u}'(x',z',\tau)|^2d\tau\,d\xi-z_3\int_{0}^1|\partial_{z_3}\widehat  {\bf u}'(x',z',\xi)|^2d\xi\right),

\end{array}
\end{equation}
where for $i=1,2$, the expression for $\widehat {\bf u}'$ is given by
\begin{equation}\label{Uhat}
\widehat {\bf u}'(x',z)={K\over \mu}\sum_{i=1}^2\left(f_i(x')-\partial_{x_i}\widetilde p(x')\right)\left(A_1^\star(z')e^{Mz_3}+A_2^\star(z')e^{-M z_3}+1\right)(e_i+\nabla_{z'}\pi^i(z')),
\end{equation}
with
\begin{equation}\label{expAZ}A_1^*(z')=-{1-e^{-Mh(z')}\over e^{Mh(z')}-e^{-Mh(z')}},\quad
A_2^*(z')={1-e^{Mh(z')}\over e^{Mh(z')}-e^{-Mh(z')}}.
\end{equation}
\end{itemize}
\end{corollary}
\begin{proof}

\noindent  We start by proving $(i)$. To prove (\ref{Expressions0ell1}), it remains to eliminate the microscopic variable $z$ in the effective problem (\ref{limit_model_0ell1}). The procedure is standard and can be viewed for instance in \cite{Allaire}, but for reader's convenience, we give some details of the proof. From the  equation (\ref{limit_model_0ell1})$_1$, the velocity $\widehat u$ is computed in terms of the macroscopic forces and the local velocities
\begin{equation}\label{identificationsub}
\widehat {\bf u}'(x',z)=-{K\over \mu}\sum_{i=1}^2\left(f_i(x')-\partial_{x_i}\widetilde p(x')\right)w^i(z),\quad \widehat \pi(x',z')=-\sum_{i=1}^2\left(f_i(x')-\partial_{x_i}\widetilde p(x')\right)\pi^i(z'),
\end{equation}
we deduce that the average velocity is given by
\begin{equation}\label{Velocity_thm_sub_proof}
\displaystyle {\bf V}_{\rm av}'(x')=-{K\over \mu}A_{M}\left({\bf f}'(x')-\nabla_{x'}\widetilde p(x')\right),\quad\hbox{in }\omega ,
\end{equation}
with $A_{M} \in \mathbb{R}^{2\times 2}$ is symmetric, positive definite and defined by its entries
\begin{equation}\label{Asub}
(A_{M})_{ij}= \int_{Z}w_j^i(z)\,dz,\quad i,j=1,2,
\end{equation}
where $w^i(z)$ $(i=1,2)$ denote the unique solution in $H^1_{\#}(Z)^2$ of the local reduced Darcy-Brinkman problem
\begin{equation}\label{localsub}
\left\{\begin{array}{ll}
\displaystyle-{1\over M^2}\partial_{z_3}^2 w^i+w^i=-(e_i+\nabla_{z'}\pi^i(z'))&\hbox{in }Z,\\
\noame
\displaystyle {\rm div}_{z'}\left(\int_0^{h(z')}w^i\,dz_3\right)=0&\hbox{in }Z',\\
\noame
\displaystyle w^i=0&\hbox{ on }\widehat \Gamma_0\cup \widehat \Gamma_1,\\
\noame
w^i,\pi^i\quad Z'-periodic.
\end{array}\right.
\end{equation}
However, we can obtain a simpler formula for $A_{M}$. We observe that local problems (\ref{localsub}) can be viewed as  systems of ordinary differential equations with constant coefficient, with respect to the variable $z_3$ and unknowns functions $z_3\to w_1^i, w_2^i$, $i=1,2$, where $z'$ is a parameter, $z'\in Z'$. The expression for $w^i$, $i=1,2$, are given by
\begin{equation}w^i(z)=-\left(A_1^\star(z')e^{Mz_3}+A_2^\star(z')e^{-M z_3}+1\right)(e_i+\nabla_{z'}\pi^i(z')),
\end{equation}
where
$$A_1^*(z')=-{1-e^{-Mh(z')}\over e^{Mh(z')}-e^{-Mh(z')}},\quad
A_2^*(z')={1-e^{Mh(z')}\over e^{Mh(z')}-e^{-Mh(z')}}.
$$
This can be computed similarly to the derivation of expressions (\ref{Expressions_thin}) from system (\ref{hom_system_sub_u_thin}) (see \cite[Corollary 2]{Pazanin_SG_Darcy} for details), so we omit the proof.

Since
$$\int_0^{h(z')}w^i(z)\,dz_3=-\left({2\over M}{e^{Mh(z')}-e^{-Mh(z')}-2\over e^{Mh(z')}-e^{-Mh(z')}}- h(z')\right)(e_i+\nabla_{z'}\pi^i(z')),\quad i=1,2,
$$
then,  from the divergence condition (\ref{localsub})$_2$, we get that $\pi^i$ is the unique solution of problem (\ref{reynoldspi}).
Therefore, the average velocity is given by (\ref{Velocity_thm_sub_proof}), with the permeability tensor $A_{M}$ defined by
$$\begin{array}{l}(A_{M})_{ij}=\displaystyle\int_{Z}w_j^i(z)\,dz= -\int_{Z'}  \left({2\over M}{e^{Mh(z')}-e^{-Mh(z')}-2\over e^{Mh(z')}-e^{-Mh(z')}}- h(z')\right)(e_i+\nabla_{z'}\pi^i(z'))e_j\,dz',
\end{array}$$
which is (\ref{AM0}). Finally, taking into account the expression for average velocity (\ref{averages0ell1})$_1$, the expression (\ref{Velocity_thm_sub_proof}) and  (\ref{relation2}), then the average velocity by (\ref{Expressions0ell1}), which concludes the proof of $(i)$.\\

We continue by proving $(ii)$.  Considering the expression of ${\bf V}_{\rm av}'$ into the divergence property (\ref{div_conv}), we prove (\ref{ReynoldsP}).
\\

We finish the proof by proving $(iii)$. We consider the equation satisfied by $\widehat T$ in (\ref{limit_model_0ell1})$_3$,  i.e.
$$\partial^2_{z_3}\widehat T=-{\mu\over   k\,K}|\widehat{\bf u}'|^2-{\mu_{\rm eff}\over k}|\partial_{z_3}\widehat  {\bf u}'|^2.$$
Integrating twice this equation with respect to $z_3$, we get
$$\widehat T(x',z)=-{\mu\over   k\,K}\int_{z_3}^1 \int_{\xi}^1|\widehat{\bf u}'(x',z',\tau)|^2\,d\tau\,d\xi-{\mu_{\rm eff}\over k}\int_{z_3}^1\int_\xi^1|\partial_{z_3}\widehat  {\bf u}'(x',z',\tau)|^2d\tau\,d\xi+B_1(x',z')z_3+B_2(x',z').$$
By using the boundary conditions $\widehat T=0$ on $\widehat \Gamma_1$, we get $B_2\equiv 0$, and using $-k\partial_{z_3}\widehat T=b$ on $\widehat \Gamma_0$, we deduce
$$B_1(x',z')=-{b\over k}-{\mu\over   kK}\int_{0}^1 |\widehat{\bf u}'(x',z',\xi)|^2\,d\xi-{\mu_{\rm eff}\over k}\int_{0}^1|\partial_{z_3}\widehat  {\bf u}'(x',z',\xi)|^2d\xi,$$
and so, expression for $\widehat T$ is given by (\ref{That}). Here, expressions for $\widehat {\bf u}'$ and $\partial_{z_3} \widehat {\bf u}'$ van be obtained by using (\ref{identificationsub}) and expressions for $w^i, i=1,2$.

Finally, taking into account the definition of the average temperature (\ref{averages0ell1})$_2$ and relation (\ref{relation2temp}), we deduce (\ref{ExpT_thm0ell1}).

\end{proof}

\subsection{The case $\ell=1$}
We start by proving that the limit functions $(\widehat{\bf u}, \widehat T, \widetilde p)$ satisfies a two-pressure problem. Then, we derive the main result for $\ell=1$ by providing the expressions for velocity $(\widetilde {\bf u}, \widetilde T)$ and the Reynolds problem satisfied by $\widetilde p$.
\begin{theorem}[Limit model]\label{thm_limit_ell1} The limit functions $\widehat{\bf  u}$, $\widehat T$ and $\widetilde  p$ given in Lemmas \ref{lem_conv_vel}, \ref{lem_conv_temp} and \ref{lem_conv_press} satisfy
\begin{equation}\label{limit_model_ell1}
\left\{\begin{array}{rl}
\displaystyle
-2\mu_{\rm eff} {\rm div}_z(\mathbb{D}_z[\widehat {\bf  u}]) +{\mu\over K}\widehat {\bf u}+\nabla_{z}\widehat \pi={\bf f}'(x') - \nabla_{x'}\widetilde  p(x')&\hbox{ in }\omega\times Z,\\
\noame
\displaystyle {\rm div}_{z'}\left(\int_0^{h(z')}\widehat {\bf u}'\,dz_3\right)=0&\hbox{ in }\omega\times Z',\\
\noame
\displaystyle-k\Delta_z\widehat T={\mu\over K}|\widehat {\bf  u}|^2+2\mu_{\rm eff}|\mathbb{D}_z[\widehat {\bf  u}]|^2 &\hbox{ in }\omega\times Z,\\
\noame

\widehat {\bf  u}'=0&\hbox{ on } \omega\times (\widehat \Gamma_0\cup \widehat \Gamma_1),\\
\noame
\widehat T=0&\hbox{ on }    \omega\times \widehat \Gamma_1,\\
\noame
k\nabla_{z}{ \widehat T}\cdot n=b&\hbox{ on }\omega\times \widehat \Gamma_0,\\
\noame
\widehat \pi\in L^{2}(\omega;L^{2}_{0,\#}(Z)).
\end{array}\right.
\end{equation}
\end{theorem}
\begin{proof}
To prove system (\ref{limit_model_ell1}), we observe that conditions (\ref{limit_model_ell1})$_{2, 4, 5}$ are consequence of Lemmas \ref{lem_conv_vel} and \ref{lem_conv_temp} in the case $\ell=1$. We divide the proof in four steps. First, we pass to the limit in the variational formulation for velocity (\ref{form_var_hat_templ1})$_1$. Then, in second and third steps, we prove strong convergences of velocity, which is used to pass to the limit in the variational formulation for temperature (\ref{form_var_hat_templ1})$_2$ in the fourth step.
\\

{\it Step 1.} Let us pass to the limit in the variational formulation for velocity (\ref{form_var_hat_templ1})$_1$:
\begin{itemize}
\item First and second terms.  From convergence of velocity given in Lemma \ref{lem_conv_vel}, we get
$$\begin{array}{l}
\displaystyle 2\mu_{\rm eff}\int_{\omega\times Z} \mathbb{D}_{z}[\widehat  {\bf u}_\ep]:\mathbb{D}_{z}[ {\bf v}]\,dx'dz\!
+\!{\mu\over  K} \int_{\omega\times Z}\!\!\widehat  {\bf u}_\ep\cdot   {\bf v}\,dx'dz\\
\noame
\displaystyle
\to 2\mu_{\rm eff}\int_{\omega\times Z} \mathbb{D}_{z}[\widehat  {\bf u}]:\mathbb{D}_{z}[ {\bf v}]\,dx'dz
\displaystyle+\!{\mu\over  K} \int_{\omega\times Z}\!\!\widehat  {\bf u}\cdot   {\bf v}\,dx'dz.
\end{array}$$
\item Third term. According to estimate (\ref{estim_p0p1})$_3$, we deduce
$$\left|\int_{\omega\times Z}\ep^2\widehat  p^\varepsilon_1 {\rm div}_{x'}({\bf v}')\,dx'dz\right|\leq C\ep\to 0,$$
and so, from convergence of pressure given in Lemma \ref{lem_conv_press}, we deduce
$$\int_{\omega\times Z}\ep^2(\widehat p_\ep^0+\widehat  p^\varepsilon_1) {\rm div}_{x'}({\bf v}')\,dx'dz\to \int_{\omega\times Z}\widetilde p\,{\rm div}_{x'}({\bf v}')\,dx'dz.$$
Observe that from condition ${\rm div}_{x'}(\int_Z {\bf v}'\,dz)=0$ in $\omega$ and since $\widetilde p$ does not depend on $z$, we get
$$\int_{\omega\times Z}\widetilde p\,{\rm div}_{x'}({\bf v}')\,dx'dz=\int_{\omega}\widetilde p\,{\rm div}_{x'}\left(\int_Z {\bf v}'dz\right)\,dx'=0.$$
\end{itemize}
Thus, from previous convergences, the limit of the variational formulation (\ref{form_var_hat_templ1})$_1$ is given by
 \begin{equation}\label{var_limit_problem_ell1}
 2\mu_{\rm eff}\int_{\omega\times Z} \mathbb{D}_{z}[\widehat  {\bf u}]:\mathbb{D}_{z}[ {\bf v}]\,dx'dz
\displaystyle+\!{\mu\over  K} \int_{\omega\times Z}\!\!\widehat  {\bf u}\cdot   {\bf v}\,dx'dz=\int_{\omega\times Z}{\bf f}'\cdot   {\bf v}'\,dx'dz,
 \end{equation}
 which, by density, holds for every ${\bf v}'\in \mathcal{V}$ with
$$\mathcal{V}=\left\{\begin{array}{l}
\displaystyle {\bf v}\in L^2(\omega;H^1_\#(Z)^3)\ :\ {\rm div}_{z}({\bf v})=0\quad\hbox{in }\omega\times Z\\
\noame
\displaystyle
{\rm div}_{x'}\left(\int_Z{\bf v}'\,dz\right)=0\quad\hbox{in }\omega,\quad \left(\int_Z{\bf v}'\,dz\right)\cdot n=0\quad\hbox{on }\partial\omega
\end{array}
\right\}.
$$
Finally, reasoning as in \cite[Lemma 1.5]{Allaire}, the orthogonal of $\mathcal{V}$ is made of gradients of the form $\nabla_{x'}\widetilde\pi(x')+\nabla_{z}\widehat \pi(x',z)$, with $\widetilde \pi\in L^{2}_0(\omega)$ and $\widehat\pi(x',z')\in L^{2}(\omega;L^{2}_{0,\#}(Z))$. Thus, integrating by parts, the limit variational formulation (\ref{var_limit_problem_ell1}) is equivalent to the two-pressures reduced Stokes system (\ref{limit_model_ell1})$_1$. It remains to prove that $\widetilde \pi$ coincides with pressure $\widetilde p$. Thus can be easily done passing to the limit similarly as above by considering the test function ${\bf v}$, which is divergence-free only in $z$, and by identifying limits.
\\

{\it Step 2.} We prove the following property
\begin{equation}\label{prevprop_crit}\begin{array}{l}
 \displaystyle
 \lim_{\ep\to 0}\left({\mu\over K}\int_{\omega\times Z}|\widehat {\bf u}_\ep|^2\,dx'dz+2\mu_{\rm eff}\int_{\omega \times Z}|  \mathbb{D}_{z}[\widehat {\bf u}_\ep]|^2 \,dx'dz\right)  ={\mu\over K}\int_{\omega\times Z}|\widehat {\bf u}|^2 \,dx'dz+2\mu_{\rm eff}\int_{\omega\times Z}|\mathbb{D}_z[\widehat {\bf u}]|^2 \,dx'dz\,.
 \end{array}
 \end{equation}
We take $\widehat {\bf u}_\ep$ as test function in  (\ref{form_var_hat_templ1})$_1$. From the fact that $\ep^{-1}{\rm div}_{z}(\widehat {\bf u}_\ep)=0$, using strong convergence of $\ep^2\widehat p_\ep^0$ and convergence of velocity,  we get
\begin{equation}\label{prevpro111crit}
\begin{array}{l}
 \displaystyle
\lim_{\ep\to 0}\left(2\mu_{\rm eff}\int_{\omega \times Z} | \mathbb{D}_{z}[\widehat {\bf u}_\ep]|^2\,dx'dz
+{\mu\over K}\int_{\omega\times Z}|\widehat {\bf u}^\ep|^2\,dx'dz\right)\\
\noame
\displaystyle=  \displaystyle\int_{\omega\times Z}\widetilde p\, {\rm div}_{x'}(\widehat {\bf u}')\,dx'dz+\int_{\omega\times Z}{\bf f}'\cdot \widehat {\bf u}'\,dx'dz=\int_{\omega\times Z}{\bf f}'\cdot \widehat {\bf u}'\,dx'dz,
\end{array}
\end{equation}
because $\int_{\omega\times Z}\widetilde p\, {\rm div}_{x'}(\widehat {\bf u}')\,dx'dz=0$, because $\widetilde p$ does not depend on $z$ and ${\rm div}_{x'}(\int_Z\widehat {\bf u}'dz)=0$.

Now, we take $\widehat {\bf u}$ as a test function in (\ref{var_limit_problem_ell1}) and deduce
\begin{equation}\label{limituhat2crit}\begin{array}{rl}
\displaystyle
\int_{\omega\times Z}{\bf f}'\cdot \widehat{\bf u}'\,dx'dz=&\displaystyle 2\mu_{\rm eff}\int_{\omega\times Z} |\mathbb{D}_z[\widehat {\bf u}]|^2\,dx'dz
+ \displaystyle {\mu\over  K}\int_{\omega\times\times Z}  |\widehat {\bf u}|^2\,dx'dz.
 \end{array}
 \end{equation}
Then, from (\ref{prevpro111crit}) and (\ref{limituhat2crit}), we get (\ref{prevprop_crit}).
\\

 {\it Step 3.} We prove the strong convergences
\begin{equation}\label{strongconvcrit}\lim_{\ep\to 0}\int_{\omega\times Z} |\widehat {\bf u}_\ep-\widehat {\bf u}|^2\,dx'dz=0,\quad \lim_{\ep\to 0}\int_{\omega\times Z} \left|\mathbb{D}_{z}[\widehat {\bf u}_\ep]-\mathbb{D}_{z}[\widehat {\bf u}]\right|^2\,dx'dz=0.
\end{equation}
To prove (\ref{strongconvcrit}), it is enough to prove that
 $$E_\ep:= 2\mu_{\rm eff}\int_{\omega\times Z} \left| \mathbb{D}_{z}[\widehat {\bf u}_\ep]-\mathbb{D}_z[\widehat {\bf u}]\right|^2\,dx'dz+{\mu\over K}\int_{\omega\times Z} |\widehat {\bf u}_\ep-\widehat {\bf u}|^2\,dx'dz\to 0.$$
 Developing the expression of $E_\ep$, we have
 $$\begin{array}{rl}
 E_\ep =&\displaystyle
 2\mu_{\rm eff}\int_{\omega\times Z} \left| \mathbb{D}_{z}[\widehat {\bf u}_\ep]\right|^2\,dx'dz+{\mu\over K}\int_{\omega\times Z} |\widehat {\bf u}_\ep|^2\,dx'dz\\
 \noame
 &\displaystyle
 +2\mu_{\rm eff}\int_{\omega\times Z} \left|\mathbb{D}_{z}[\widehat {\bf u}]\right|^2\,dx'dz-4\mu_{\rm eff}\int_{\omega\times Z} \mathbb{D}_{z}[\widehat {\bf u}_\ep]:\mathbb{D}_{z}[\widehat {\bf u}] \,dx'dz\\
 \noame
 &\displaystyle
 +{\mu\over K}\int_{\omega\times Z} |\widehat {\bf u}|^2\,dx'dz-2{\mu\over K}\int_{\omega\times Z} \widehat {\bf u}_\ep\cdot \widehat {\bf u}\,dx'dz.
 \end{array}$$
 By using property (\ref{prevprop_crit})  and convergences of velocity, we get
 $$\begin{array}{rl}
 E_\ep \to &\displaystyle
 2\mu_{\rm eff}\int_{\omega\times Z} \left|\mathbb{D}_z[\widehat {\bf u}]\right|^2\,dx'dz+{\mu\over K}\int_{\omega\times Z} |\widehat {\bf u}|^2\,dx'dz\\
 \noame
 &\displaystyle
 +2\mu_{\rm eff}\int_{\omega\times Z} \left|\mathbb{D}_z[\widehat {\bf u}]\right|^2\,dx'dz-4\mu_{\rm eff}\int_{\omega\times Z} |\mathbb{D}_z[\widehat {\bf u}]|^2\,dx'dz\\
 \noame
 &\displaystyle
 +{\mu\over K}\int_{\omega\times Z} |\widehat {\bf u}|^2\,dx'dz-2{\mu\over K}\int_{\omega\times Z} |\widehat {\bf u}|^2\,dx'dz=0,
 \end{array}$$
 which proves the strong convergences (\ref{strongconvcrit}). \\

{\it Step 4.} Let us pass to the limit in  the variational formulation for temperature (\ref{form_var_hat_templ1})$_2$:
\begin{itemize}
\item First  term. From convergence of temperature given in Lemma  \ref{lem_conv_temp}, we get
 $$\int_{\omega\times Z} k\nabla_{z}\widehat  T_\ep\cdot \nabla_{z}  \zeta\,dx'dz
 \to \int_{\omega\times Z}
k\nabla_z\widehat T\cdot\nabla_z\zeta\,dx'dz.$$

\item Second and third terms.  From the strong convergences (\ref{strongconvcrit}), we deduce
$${\mu\over K} \int_{\omega\times Z} |\widehat  {\bf u}_\ep|^2  \zeta\,dx'dz+2\mu_{\rm eff} \int_{\omega\times Z}|\mathbb{D}_{z}[\widehat  {\bf u}_\ep]|^2  \zeta\,dx'dz\to {\mu\over K} \int_{\omega\times Z} |\widehat  {\bf u}|^2  \zeta\,dx'dz+2\mu_{\rm eff} \int_{\omega\times Z}|\mathbb{D}_{z}[\widehat  {\bf u}]|^2  \zeta\,dx'dz.$$
\end{itemize}
Thus, from previous convergences,  the limit of the variational formulation (\ref{form_var_hat_subcritical})$_2$ is given by
 \begin{equation}\label{var_limit_problem_temp_0ell1}\int_{\omega\times Z}
k\nabla_{z}\widehat T\cdot\nabla_{z}\zeta\,dx'dz=2\mu_{\rm eff}\int_{\omega\times Z}|D_z[\widehat{\bf u}]|^2\,dx'dz
\displaystyle+\!{\mu\over  K} \int_{\omega\times Z}|\widehat  {\bf u}|^2\,dx'dz+\int_{\omega\times Z'}   b\,  \zeta\,dx'dz',
\end{equation}
which, by density, holds for every $\zeta\in L^q(\omega;W^{1,q}_\#(Z))$. Thus, integrating by parts, the limit variational formulation (\ref{var_limit_problem_temp_0ell1}) is equivalent to the two-pressures reduced Stokes system (\ref{limit_model_ell1})$_{3,6}$.

\end{proof}

\begin{corollary} Consider $(\widehat u, \widetilde p, \widehat T)$, with $\widehat u_3=0$, satisfying system (\ref{limit_model_0ell1}). Let us define the average velocity and temperature by
\begin{equation}\label{averages0ell1}
{\bf V}_{\rm av}(x')=\int_0^{h_{\rm max}}\widetilde {\bf u}(x',z_3)\,dz_3,\quad T_{\rm av}(x')=\int_0^{h_{\rm max}}\widetilde T(x',z_3)\,dz_3,
\end{equation}
and a parameter $M=\sqrt{\mu\over K\mu_{\rm eff}}>0$.  We have the following:
\begin{itemize}
\item[(i)] The average velocity ${\bf V}_{\rm av}$ is given by
 \begin{equation}\label{Velocity_thm_crit}
\begin{array}{lll}
\displaystyle {\bf V}_{av}'(x')={1\over M^2 \mu_{\rm eff}}A_{M}\left({\bf f}'(x')-\nabla_{x'}\widetilde p(x')\right),\quad&  {V}_3^{av}(x')=0,\quad&\hbox{in }\omega,
\end{array}
\end{equation}
where   $A_{M} \in \mathbb{R}^{2\times 2}$, which is symmetric and positive definite, is defined by its entries
\begin{equation}\label{Acrit}
(A_{M})_{ij}= \int_{Z}\mathbb{D}_\lambda[w^i]: \mathbb{D}_\lambda[w^j]\,dz,\quad i,j=1,2,
\end{equation}
with $w^i(z),  i=1,2$  the unique solution in $H^1_{\#}(Z)^3$ of the local 3D Darcy-Brinkman problem
\begin{equation}\label{localcrit}
\left\{\begin{array}{ll}
\displaystyle-{1\over M^2}{\rm div}(\mathbb{D}[w^i])+\nabla\pi^i+w^i=e_i&\hbox{in }Z,\\
\noame
{\rm div}_z(w^i)=0&\hbox{in }Z,\\
\noame
\displaystyle w^i=0&\hbox{ on }\widehat \Gamma_0\cup \widehat \Gamma_1,\\
\noame
w^i,\pi^i\quad Z'-periodic.
\end{array}\right.
\end{equation}
\item[(ii)] Pressure $\widetilde p \in L^2_0(\omega)\cap H^1(\omega)$ is the unique solution of the Reynolds equation:
\begin{equation}\label{Pcrit}
 \left\{\begin{array}{rl}
\displaystyle {\rm div}_{x'}\left({1\over M^2 \mu_{\rm eff}}{A_M}\left({\bf f}'(x')-\nabla_{x'}\widetilde p(x')\right)\right)=0& \hbox{in } \omega,
\\
\noame
\displaystyle\left({1\over M^2 \mu_{\rm eff}}{A_M}\left({\bf f}'(x')-\nabla_{x'}\widetilde p(x')\right)\right)\cdot n=0& \hbox{on } \partial\omega.
\end{array} \right.\end{equation}

\item[(iii)] The average temperature $T_{\rm av}$ is given by
$$T_{av}(x')=\int_Z\widehat T(x',z)\,dz\quad \hbox{in }\omega ,$$
with $\widehat T\in L^q(\omega ;W^{1,q}_{\#}(Z))$ the unique solution of the nonlinear problem
\begin{equation}
\left\{\begin{array}{rl}\label{identificationaa}\displaystyle
-k\Delta_z \widehat T={\mu\over K}\left|\widehat {\bf u}\right|^2+2\mu_{\rm eff}\left|\mathbb{D}_z[\widehat {\bf u}]\right|^2&\hbox{in }\omega\times Z,\\
\noame
\widehat T=0&\hbox{ on }    \omega\times \widehat \Gamma_1,\\
\noame
k\nabla_{z}{ \widehat T}\cdot n=b&\hbox{ on }\omega\times \widehat \Gamma_0,
\end{array}\right.
\end{equation}
where
\begin{equation}\label{identification2}
\widehat {\bf u}(x',z)={K\over \mu}\sum_{i=1}^2\left(f_i(x')-\partial_{x_i}\widetilde p(x')\right)w^i(z).
\end{equation}

\end{itemize}
\end{corollary}

\begin{proof}
We start b proving $(i)$. To prove (\ref{Velocity_thm_crit}), it remains to eliminate the microscopic variable $z$ in the effective problem (\ref{limit_model_ell1}). The procedure is standard and can be viewed for instance in \cite{Allaire}, but for reader's convenience, we give some details of the proof. From the  equation (\ref{limit_model_ell1})$_1$, the velocity $\widehat {\bf u}$ is computed in terms of the macroscopic forces and the local velocities
\begin{equation}\label{identification}
\widehat {\bf u}(x',z)={K\over \mu}\sum_{i=1}^2\left(f_i(x')-\partial_{x_i}\widetilde p(x')\right)w^i(z),\quad \widehat \pi(x',z)=\sum_{i=1}^2\left(f_i(x')-\partial_{x_i}\widetilde p(x')\right)\pi^i(z),
\end{equation}
with $(w^i, \pi^i)$, $i=1,2$, solving the local problem (\ref{localcrit}). \\

Integrating this expression on $Z$ and taking into account that ${\bf V}_{\rm av}(x')=\int_0^{h_{\rm max}}\widetilde{\bf u}(x',z_3)\,dz_3=\int_{Z}\widehat{\bf u}(x',z)\,dz$ and $\int_Z\widehat u_3\,dz=0$, we get the relationship (\ref{Velocity_thm_crit}), because   the matrix $A_{M}\in\mathbb{R}^{2\times 2}$ satisfies
$$(A_{M})_{ij}=\int_Z w^i_j(z)\,dz=\int_Z \mathbb{D}[w^i(z)]:\mathbb{D}[w^j(z)]\,dz,\quad i,j=1,2.$$

We continue by proving $(ii)$. Considering the expression of ${\bf V}_{\rm av}$ into the divergence property ${\rm div}_{x'}(\int_0^{h_{\rm max}}\widetilde {\bf u}'\,dz_3)=0$ in $\omega$ and $(\int_0^{h_{\rm max}}\widetilde {\bf u}'\,dz_3)cdot n=0$ on $\partial\omega$,  we prove (\ref{Pcrit}).
\\

We finish the proof by proving $(iii)$. From the expression of $\widehat {\bf u}$ given in (\ref{identification}) and problem satisfied by $\widehat T$ given in (\ref{limit_model_ell1}), we deduce that $T_{\rm av}$ satisfies problem (\ref{identificationaa})-(\ref{identification2}).

\end{proof}

\subsection{The case $\ell>1$}
In this section, we comment on the case of highly oscillating boundary, i.e.~ for $\ell>1$ where the thickness of the domain $\ep$ is greater than the period $\ep^\ell$. We introduce the following notation:
\begin{itemize}
\item  $\widetilde \Omega_+^\ep=\omega\times (h_{\rm min}, h(x'/\ep^\ell))$ the oscillating part of the domain,
\item $\Omega_+=\omega\times (h_{\rm min}, h_{\rm max})$ the extension of the oscillating part,
\item $\widetilde {\bf u}_\ep^+$ and $\widetilde T_\ep^+$ the restriction  to $\widetilde \Omega_+^\ep$ of the velocity and temperature. We denote by the same symbols the extension to $\Omega_+$, because the homogeneous boundary conditions on the top  boundary.
\end{itemize}
\paragraph{Estimates for velocity and temperature.} Because the following Poincar\'e's inequality in $\Omega^\ep_+$ given in \cite[Lemma 4.3]{Anguiano_SG}:
$$\|w\|_{L^r(\Omega^\ep_+)}\leq C\ep^\ell\|\nabla_{x'}w\|_{L^r(\Omega^\ep_+)^2},\quad \forall\,w\in W^{1,r}(\Omega^\ep_+),\ w_{|_{\widetilde\Gamma^\ep_1}}=0,\quad 1\leq r<+\infty,$$
it can be proved that the following estimates hold:
$$\|\widetilde {\bf u}_\ep^+\|_{L^2(\widetilde\Omega_+^\ep)^3}\leq C\ep^{2(\ell-1)},\quad\|D_\ep\widetilde {\bf u}_\ep^+\|_{L^2(\widetilde\Omega_+^\ep)^3}\leq C\ep^{\ell-2}.$$
 Moreover, following the derivation of the estimates for temperature in \cite[Proposition 1]{Pazanin_SG_Darcy} applied for $\widetilde T_\ep^+$ in $\widetilde \Omega_+^\ep$ and using the Poincar\'e inequality, we can derive the following estimates
$$\|\widetilde T_\ep\|_{L^q(\widetilde\Omega^\ep_+)}\leq C\ep^{\ell-1+{2\over q}(\ell-1)(3-q)},\quad \|\nabla_\ep \widetilde T_\ep\|_{L^q(\widetilde\Omega^\ep_+)^3}\leq C\ep^{-1+{2\over q}(\ell-1)(3-q)}.$$
The same estimates also hold for the extensions to $\Omega_+$.
\paragraph{Convergences for velocity and temperature.} From the previous estimates for the extended $\widetilde {\bf u}^+_\ep$ and $\widetilde T_\ep^+$ to $\Omega_+$, we deduce the following convergences
\begin{equation}\label{conv_high}\widetilde {\bf u}_\ep^+\rightharpoonup 0\quad\hbox{in }L^2(\Omega_+)^3,\quad \widetilde T_\ep^+\rightharpoonup 0\quad\hbox{in }L^q(\Omega_+).
\end{equation}
The proof follows from the equality
$$\widetilde {\bf u}_\ep^+=\ep^{2(\ell-1)}\left(\ep^{-2(\ell-1)}\widetilde {\bf u}_\ep^+\right),$$
and so, since $\ep^{-2(\ell-1)}\widetilde {\bf u}_\ep^+$ is bounded and $\ep^{2(\ell-1)}\to 0$, because $\ell>1$, we get the convergence for velocity.  Similarly,
$$\widetilde T_\ep^+=\ep^{\ell-1+{2\over q}(\ell-1)(3-q)}\left(\ep^{-(\ell-1+{2\over q}(\ell-1)(3-q))}\widetilde T_\ep^+\right),$$
and so, since $\ep^{-(\ell-1+{2\over q}(\ell-1)(3-q))}\widetilde T_\ep^+$ is bounded and $\ep^{\ell-1+{2\over q}(\ell-1)(3-q)}\to 0$, because $\ell-1+{2\over q}(\ell-1)(3-q)>0$, we get the convergence for temperature.

The case $\ell>1$, due to the highly oscillating boundary, leads to the conclusion that
the velocity and temperature are zero in the oscillating part $\Omega_+$. So, according to this, at macroscopic level, our problem reduces to the study of the fluid flow in the thin domain $\Omega_{\rm min}=\omega\times (0,h_{\rm min})$, which is, in fact, a thin domain without roughness considered in \cite{Pazanin_SG_Darcy}. For reader's convenience, we provide below the expressions for the velocity and temperature and the corresponding Reynolds equation for pressure obtained in \cite{Pazanin_SG_Darcy}, just replacing $h(x')$ by $h_{\rm min}$.

Denoting  the $\ep$-independent domain $\Omega=\omega\times (0,h(x'))$ with bottom boundary $\Gamma_0=\omega\times\{0\}$ and top boundary $\Gamma_1=\omega\times \{h(x')\}$, it was deduced the following convergence result given in \cite[Proposition 3]{Pazanin_SG_Darcy}:
\begin{itemize}
\item There exists  $\widetilde {\bf u}=(\widetilde{\bf  u}^\star,0)$ with $\widetilde{\bf  u}^\star\in V^2_{z_3}(\Omega)^2$ and $\widetilde{\bf  u}^\star=0$ on $\Gamma_0\cup \Gamma_1$ such that
$$\begin{array}{c}
\displaystyle \widetilde {\bf u}_\ep\rightharpoonup \widetilde{\bf  u}^\star\quad \hbox{in } V_{z_3}^2(\Omega)^2,\\
\noame
\displaystyle
 \widetilde {u}_{\ep,3}\rightharpoonup0\quad \hbox{in } L^2(\Omega),\quad
 \ep \partial_{x_j}\widetilde {\bf u}'_\ep\rightharpoonup 0 \quad\hbox{in }L^2(\Omega),\quad j=1,2,\\
 \noame
 \displaystyle
 {\rm div}_{x'}\left(\int_0^{h(x')}\widetilde{\bf  u}^\star(x',z_3)\,dz_3\right)=0\quad\hbox{in }\omega,
 \quad
 \left(\int_0^{h(x')}\widetilde{\bf  u}^\star(x',z_3)\,dz_3\right)\cdot n=0\quad\hbox{on }\partial\omega.
 \end{array}
$$
\item There exists $\widetilde p^\star \in L^2_0(\omega)\cap H^1(\omega)$ such that
$$\ep^2\widetilde p_\ep^0\rightharpoonup \widetilde p^\star\quad\hbox{in }H^1(\omega).$$
\item There exists $\widetilde T^\star \in V_{z_3}^q(\Omega)$, $q\in (1,3/2)$ such that
$$\widetilde T_\ep\rightharpoonup \widetilde T^\star\quad\hbox{in }V_{z_3}^q(\Omega),\quad \ep \partial_{x_j}\widetilde T_\ep\rightharpoonup 0\quad\hbox{in }L^q(\Omega),\ j=1,2.$$

\end{itemize}
Moreover, the limit functions $\widetilde{\bf  u}^\star$, $\widetilde  p^\star$ and  $\widetilde T^\star$  satisfy the reduced system (see \cite[Theorem 1]{Pazanin_SG_Darcy} for more details):
\begin{equation}\label{hom_system_sub_u_thin}
\left\{\begin{array}{rl}
\displaystyle
-\mu_{\rm eff} \partial_{z_3}^2 \widetilde {\bf  u}^\star +{\mu\over K}\widetilde {\bf u}^\star={\bf f}'(x') - \nabla_{x'}\widetilde  p^\star(x')&\hbox{ in }\Omega,\\
\displaystyle {\rm div}_{x'}\left(\int_0^{h(x')}\widehat {\bf u}^\star\,dz_3\right)=0&\hbox{ in }\Omega,\\
\displaystyle-k\partial^2_{z_3}\widetilde T^\star={\mu\over K}|\widetilde {\bf  u}^\star|^2+\mu_{\rm eff}|\partial_{z_3}\widetilde {\bf  u}^\star|^2 &\hbox{ in }\Omega,\\
\widetilde {\bf  u}^\star=0&\hbox{ on } \Gamma_0\cup \Gamma_1,\\
\widetilde T^\star=0&\hbox{ on }    \Gamma_1,\\
-k\partial_{z_3}{ \widetilde T^\star}=b&\hbox{ on }\Gamma_0,
\end{array}\right.
\end{equation}
and the following expressions for $\widetilde {\bf u}^\star$, $\widetilde p^\star$ and $\widetilde T^\star$ hold (see \cite[Corollary 2]{Pazanin_SG_Darcy} for more details):
\begin{itemize}
\item  Velocity $\widetilde {\bf u}^\star$ is given by
 \begin{equation}\label{Expressions_thin}
\widetilde {\bf u}^\star(x',z_3)={K\over \mu}\left(A_1^\star(x')e^{Mz_3}+A_2^\star(x')e^{-M z_3}+1\right)({\bf f}'(x')-\nabla_{x'}\widetilde p^\star(x')),
\end{equation}
where $M=\sqrt{\mu\over K\mu_{\rm eff}}$ and
 \begin{equation}\label{expA_thin}A_1^*(x')=-{1-e^{-Mh(x')}\over e^{Mh(x')}-e^{-Mh(x')}},\quad
A_2^*(x')={1-e^{Mh(x')}\over e^{Mh(x')}-e^{-Mh(x')}}.
\end{equation}
\item Pressure $\widetilde p^\star$ is the unique solution of the Reynolds equation:
\begin{equation}\label{ReynoldsP_thin}
\left\{ \begin{array}{l}\displaystyle {\rm div}_{x'}
\left({K\over \mu}\left({2\over M}{e^{Mh(x')}-e^{-Mh(x')}-2\over e^{Mh(x')}-e^{-Mh(x')}}- h(x')\right)(\nabla_{x'}\widetilde p^\star-{\bf f}')\right)=0\quad\hbox{in }\omega.\\
\noame
\displaystyle
\left({K\over \mu}\left({2\over M}{e^{Mh(x')}-e^{-Mh(x')}-2\over e^{Mh(x')}-e^{-Mh(x')}}-  h(x')\right)(\nabla_{x'}\widetilde  p^\star-{\bf f}')\right)\cdot n=0\quad\hbox{on }\partial\omega.
\end{array}\right.
\end{equation}

\item The temperature $\widetilde T^\star$ is given by
\begin{equation}\label{ExpT_thm_thin}
\begin{array}{rl}\displaystyle \widetilde T^\star(x',z_3)=&-\displaystyle {K\over k\mu}\Big( V_1^\star(x',z_3)-V_1^\star(x',h(x'))\Big)|{\bf f}'(x')-\nabla_{x'}\widetilde p^\star(x')|^2\\
\noame
&\displaystyle  -{\mu_{\rm eff} K^2 M^2\over k \mu^2}\Big(V_2^\star(x',z_3)-V_2^\star(x',h(x'))\Big)|{\bf f}'(x')-\nabla_{x'}\widetilde p^\star(x')|^2\\
\noame
&\displaystyle-{b\over k}(z_3-h(x'))|{\bf f}'(x')-\nabla_{x'}\widetilde p^\star(x')|^2,
\end{array}\end{equation}
with
\begin{equation}\label{V1star_thin}
\begin{array}{rl}
\displaystyle
 \medskip
V_1^\star(x',z_3)
=&\displaystyle {1\over 4M^2}\Big(A_1^\star(x')^2\left(e^{2Mz_3}-1\right)+A_2^\star(x')^2\left(e^{-2Mz_3}-1\right)\Big)\\
 \medskip
&\displaystyle+{2\over M^2}\Big(A_1^\star(x')(e^{Mz_3}-1)+A_2^\star(x')(e^{Mz_3}-1)\Big)\\
  \medskip
&\displaystyle + \left({1\over 2}+A_1^\star (x')A_2^\star(x')\right)z^2_3-{1\over 2M}\Big(A_1^\star(x')^2-A_2^\star(x')^2\Big)z_3\\
 \medskip
&\displaystyle-{2\over M}(A_1^\star(x')-A_2^\star(x'))z_3,
\end{array}
\end{equation}
\begin{equation}\label{V2star_thin}\begin{array}{rl}
\displaystyle\medskip
V_2^\star(x',z_3)= &\displaystyle{1\over 2}\left((A_1^\star)^2\left(e^{2Mh(x')}-e^{2Mz_3}\right)-(A_2^\star)^2\left(e^{-2Mh(x')}-e^{-2Mz_3}\right)\right) \\
\medskip
&\displaystyle+A_1^\star A_2^\star (h(x')-z_3)^2.
\end{array}
\end{equation}
\end{itemize}

\subsection*{Acknowledgements}
The first author  belongs to the ``Mathematical Analysis" Research Group (FQM104) at Universidad de Sevilla.
The second author is supported by the Croatian Science Foundation under the project AsyAn (IP-2022-10-1091) and Republic of Croatia's MSEY in course of Multilateral scientific and technological cooperation in Danube region under the project MultiHeFlo. \ \\
\ \\

\end{document}